\documentclass[12pt,draft]{extractic}

\usepackage{amsmath,amsthm,amssymb,calc}
\usepackage[dvips]{graphics, graphicx}

\usepackage[english]{babel}

\def\eps{\varepsilon}

\newcounter{num}[section]

\newcommand{\Th}{\refstepcounter{num}
{\bf Theorem \arabic{section}.\arabic{num} }}

\newcommand{\Lemma}{\refstepcounter{num}
{\bf Lemma \arabic{section}.\arabic{num} }}

\newcommand{\Pred}{\refstepcounter{num}
{\bf Proposition \arabic{section}.\arabic{num} }}

\newcommand{\Cor}{\refstepcounter{num}
{\bf Corollary \arabic{section}.\arabic{num} }}

\newcommand{\Note}{\refstepcounter{num}
{\it Note \arabic{section}.\arabic{num} }}

\newcommand{\Def}{\refstepcounter{num}
{\it Definition \arabic{section}.\arabic{num} }}

\newcommand{\Proof}{{\bf Proof. }}

\def\eps{\varepsilon}
\def\_phi{\varphi}

\def\a{\alpha}
\def\d{\delta}
\def\l{\lambda}

\def\v{\vec}
\def\F{\widehat}
\def\L{\Lambda}
\def\m{\times}
\def\t{\tilde}
\def\o{\omega}
\def\ov{\overline}

\def\C{{\mathbb C}}
\def\r{\mathcal{R}}
\def\Z_N{{\mathbb Z}_N}
\def\Z{{\mathbb Z}}

\def\per{{\rm per\,}}
\def\supp{{\rm supp\,}}

\author{Shkredov I.D.}
\title{On sumsets of dissociated sets
\footnote{This work was supported by National Science Foundation grant No. DMS--0635607.}}
\date{}
\begin{document}
\maketitle

\begin{center}
 Annotation.
\end{center}

{\it \small
    In the paper we are studying some properties of subsets $Q$ of sums of dissociated sets.
    The exact upper bound for the number of solutions of the following equation
    \begin{equation}\label{f:annot_1}
        q_1 + \dots + q_p = q_{p+1} + \dots + q_{2p} \,, \quad q_i \in Q
    \end{equation}
    in groups $\mathbf{F}_2^n$ is found.
    Using our approach, we easily prove a recent result of J. Bourgain
    on sets of large exponential sums and obtain a tiny improvement of his theorem.
    Besides an inverse problem is considered in the article.
    Let $Q$ be a set belonging a sumset of two dissociated sets
    such that equation (\ref{f:annot_1}) has many solutions.
    We prove that in the case the large proportion of $Q$ is highly structured.
}
\\
\\
\\

\refstepcounter{section}
\label{introduction}

\begin{center}
{\large
{\bf \arabic{section}. Introduction.}}
\end{center}

Let $G=(G,+)$ be a finite Abelian group with additive group operation $+$.
Suppose that $A$ is a subset of $G$.
It is very convenient to write $A(x)$ for such a function.
Thus $A(x)=1$ if $x\in A$ and $A(x)=0$ otherwise.
By $\F{G}$ denote the Pontryagin dual of $G$,
in other words the space of homomorphisms $\xi$ from $G$ to
$\mathbf{R} / \mathbf{Z}$,
$\xi : x \to \xi \cdot x$.
It is well known that $\F{G}$ is an additive group
which is isomorphic to $G$.
Also denote by $N$ the cardinality of $G$.
Let $f$ be a function from $G$ to  $\mathbf{C}$, $N = |G|$.
By  $\F{f}(\xi)$ denote the Fourier transformation of $f$
\begin{equation}\label{}
  \F{f}(\xi) =  \sum_{x \in G} f(x) e( -\xi \cdot x) \,,
\end{equation}
where $e(x) = e^{2\pi i x}$ and $\xi \in \F{G}$.

Let $\d,\a$ be real numbers, $0<\a \le \d \le 1$ and
let $A$ be a subset of $\Z_N$ of cardinality $\d N$.
Consider the set $\r_\a$ of large exponential sums of the set $A$
\begin{equation}\label{f:R_def}
    \r_\a = \r_\a (A) = \{~ r\in \F{G} ~:~ |\F{A} (r)| \ge \a N ~\} \,.
\end{equation}
In many problems of combinatorial number theory is important to know the structure of the set $\r_\a$
(see \cite{Gow_surv}).
In other words what kind of properties $\r_\a$ has?
Clearly,  that this question is an inverse problem
of additive number theory (see \cite{Nathanson,Tao_Vu_book}).

The first non--trivial result in the direction
was obtained by M.--\,C. Chang \cite{Ch_Fr} in 2002.
Recall that a set $\mathcal{D} = \{ d_1, \dots, d_{|\mathcal{D}|} \} \subseteq G$
is called {\it dissociated} if any equality of the form
\begin{equation}\label{f:def_diss}
        \sum_{i=1}^{|\mathcal{D}|} \eps_i d_i = 0 \,,
\end{equation}
where $\eps_i \in \{ -1,0,1 \}$ implies that all $\eps_i$ are equal to zero.

Let $\log$ stand for the logarithm to to base $2$.
Let  $p$ be a positive integer.
By $[p]$ denote the segment of natural numbers $\{1,\dots,p\}$.

 \Th {\bf (Chang)}
 \label{t:Chang}
 { \it
    Let $\d,\a$ be real numbers, $0<\a \le \d \le 1$,
    $A$ be a subset of $G$, $|A|= \d N$,
    and the set $\r_\a$ is defined by  (\ref{f:R_def}).
    Then any dissociated set $\L$, $\L \subseteq \r_\a$
    has the cardinality at most $2 ( \d/\a )^2 \log (1/\d)$.
 }

 A simple consequence of Parseval's identity gives $|\L| \le \d/\a^2$.
 Hence Chang's Theorem is nontrivial if $\d$ is small.

Using approach of paper \cite{Ruzsa_Freiman} (see also  \cite{Bilu})
 Chang applied her result to prove the famous Freiman's Theorem \cite{Freiman}
 on sets with small doubling.
 Another applications of Theorem \ref{t:Chang} were obtained by B. Green in paper \cite{GreenA+A}
 by B. Green and I. Ruzsa in \cite{Ruzsa_Green_Fr},
 T. Sanders (see e.g.  \cite{Sanders_appl_Chang,Sanders_AP3_sumsets,Sanders_Littlewood-Gowers}),
 and also T. Schoen in \cite{Schoen}.
  If the parameter $\a$ is close to $\d$ then the structural properties of the set $\r_\a$
 were studied in papers
 \cite{Freiman_Yudin,Besser,Lev} (see also survey \cite{Konyagin_Lev}).

 By $A_1 \dotplus A_2 \dotplus \dots \dotplus A_d$
 denotes the set of sums of different elements of the sets $A_1,\dots,A_d$.
 If all $A_i$ are equal to $A$ then we shall write $d\dot{A}$.

 In paper \cite{Bourgain_AP2007} J. Bourgain used sumsets of a dissociated set $\L$
 and obtained an extension
 of Chang's theorem.
 He used the extension in proving of his beautiful result on density
 of subsets of $[N]$ without arithmetic progressions of length three.
 Further applications on the Theorem below were obtained in \cite{Sanders_Fr_new}.

 \Th {\bf (Bourgain)}
 \label{t:Bougrain_new}
 { \it
    Let $d$ be a positive integer,
    $\d,\a$ be real numbers, $0<\a \le \d \le 1$,
    $A$ be a subset of $G$, $|A| = \d N$,
    and the set $\r_\a$ is defined by  (\ref{f:R_def}).
    Suppose that $\L$ is a dissociated set.
    Then for any $d\ge 1$, we have
    $|d\dot{\L} \bigcap \r_\a| \le 8 ( \d/\a )^2 \log^d (1/\d)$.
 }

 In articles \cite{Sh_dokl_exp1,Sh_exp1,Sh_exp2} another results on
 sets of large exponential sums were obtained.
 In particular, the following theorem was proved in these papers.

 \Th
 \label{t:main}
 {\it
    Let $\d,\a$ be real numbers, $0< \a \le \d$,
    $A$ be a subset of $\Z_N$, $|A| = \d N$, and
    $k\ge 2$ be a positive integer.
    Let also $B\subseteq \r_\a \setminus \{ 0 \}$ be an arbitrary set.
    Then the number
    \begin{equation}\label{f:T_k_def}
        T_k (B) := |\{ ~ (r_1,\dots, r_k, r_1', \dots, r_k') \in B^{2k} ~:~
                    r_1 + \dots + r_k = r_1' + \dots + r_k' ~ \}|
    \end{equation}
    is at least
    \begin{equation}\label{f:T_k}
        \frac{\d \a^{2k}}{2^{4k} \d^{2k}} |B|^{2k} \,.
    \end{equation}
 }

 In article \cite{Sh_exp1} was showed that
 Theorem \ref{t:main} and an inequality of W. Rudin \cite{Rudin_book,Rudin}
 on dissociated sets
 imply M.--\,C. Chang's theorem.
 Similarly in the paper we show that an appropriate analog
 of Rudin's result and Theorem \ref{t:main} gives us Theorem \ref{t:Bougrain_new}
 in $\mathbf{F}_2^n$
 (see section \ref{elem_Bourgain}).
 Our approach is an elementary and does not require
 sufficiently difficult hypercontractivity technic from  \cite{Bourgain_AP2007}.
 We show that for any $Q \subseteq d \dot{\L}$,
 where $\L$ is a dissociated,
 the value $T_k (Q)$ does not exceed $C^{dk} k^{dk} |Q|^{k}$.
 Here
 $C>0$ is an absolute constant.
 Applying this result to the set $d\dot{\L} \bigcap \r_\a$
 and using Theorem \ref{t:main}, we get Theorem \ref{t:Bougrain_new}.
 Actually a tiny improvement of the last result was obtained
 (see Theorem \ref{t:my_Bougrain_new}).

 In section \ref{inverse} an inverse problem is considered.
 Let $Q$ be a subset of $2 \dot{\L}$, where
 $\L$ is an arbitrary dissociated set.
 Suppose that the value $T_k (Q)$ is large
 in the sence that $T_k(Q) \gg C^{dk} k^{dk} |Q|^{k}$.
 What can we say about the structure of $Q$?
 We show that in the case the set $Q$ contains a sumset of two dissociated sets
 (see Theorem \ref{t:inverse}).
 In some sence we give a full description of large subsets of
 $2\dot{\L}$ with large value of $T_k$.

 The obtained results are formulated in groups $\mathbf{F}^n_2$
 but they can be extended to any Abelian group
 (see discussion of using $\mathbf{F}^n_p$, $p$ is a prime,
  in \cite{Green_finite_fields}).
In our forthcoming papers we are going to obtain these extensions.


I
acknowledge the Institute for Advanced Study for its hospitality
and providing me with excellent working conditions.
Also the author is grateful to Professor N.G. Moshchevitin and Professor S.V. Konygin
for attention to this work and useful discussions.

\refstepcounter{section}
\label{elem_Bourgain}

\begin{center}
{\large
{\bf \arabic{section}. An elementary proof of a result of Bourgain.}}
\end{center}

Denote by $G$ the group $\mathbf{F}_2^n$.
Let $A\subseteq G$ be a set, and $k\ge 2$ be a positive integer.
By $T_k (A)$ denote the number
$$
    T_k (A) := | \{ a_1 + \dots + a_k = a'_1 + \dots + a'_k  ~:~ a_1, \dots, a_k, a'_1,\dots,a'_k \in A \} | \,.
$$
If $A_1,\dots,A_{2k} \subseteq G$ are any sets,
then denote by $T_k (A_1,\dots,A_{2k})$ the following number
$$
    T_k (A_1,\dots,A_{2k}) := | \{ a_1 + \dots + a_k = a_{k+1} + \dots + a_{2k}  ~:~ a_i \in A_i,\, i=1,\dots, 2k \} | \,.
$$
We shall write $\sum_x$ instead of $\sum_{x\in G}$ for simplicity.

    Using the notion of convolution,
    we can calculate $T_k (A)$.

\Def
\label{def:convolution}
    Let $f,g : G \to \mathbb{C}$ be any  functions.
    Denote by $(f*g) (x)$ the function
\begin{equation}\label{f:*-conv}
    (f*g) (x) = \sum_s f(s) g(x-s) \,.
\end{equation}
    Clearly,  $(f*g) (x) = (g*f) (x)$, $x\in G$.
    Further, using induction, we get the operation $*_k$, where $k$ is a positive integer.
    So $*_k := *(*_{k-1})$.

If $A,B \subseteq G$ are arbitrary sets, then $(A*B) (x) \neq 0$
iff  $x\in A+B$.
Hence $T_2 (A) = \sum_x (A*A)^2 (x)$.
Let $f: G \to \mathbb{C}$ be a function.
By $T_k (f)$ denote
$T_k (f) = \sum_x |(f *_{k-1} f) (x)|^2$.

\Lemma
\label{l:T_k(B,A)}
{\it
    Let $s$,$t$ be positive integers, $s\ge 2$, $t\ge 2$, and let
    $f_1, \dots, f_{s}, g_1, \dots, g_{t} : G \to \mathbb{R}$ be functions.
    Then
    $$
        \left| \sum_x (f_1 * \dots * f_{s}) (x) \cdot (g_1 * \dots * g_{t}) (x) \right|
            \le
    $$
    \begin{equation}\label{f:T_k(B,A)}
            \le
                (T_{s} (f_1))^{1/2s} \dots (T_{s} (f_{s}))^{1/2s}
                    (T_{t} (g_1))^{1/2t} \dots (T_{t} (g_{t}))^{1/2t} \,.
    \end{equation}
}
\label{l:conv}
\Proof
Since $\F{(f*g)} (r) = \F{f} (r) \F{g} (r)$, it follows that
$$
    \sigma := \sum_x (f_1 * \dots * f_{s}) (x) \cdot (g_1 * \dots * g_{t}) (x)
        =
            \frac{1}{N} \sum_r \F{f}_1 (r) \dots \F{f}_{s} (r) \ov{\F{g}_1 (r)} \dots \ov{\F{g}_{t} (r)} \,.
$$
Using H\"{o}lder's inequality several times, we obtain
$$
    \sigma
        \le
            \left( \frac{1}{N} \sum_r |\F{f}_1 (r)|^{2s} \right)^{\frac{1}{2s}}
                \dots
                    \left( \frac{1}{N} \sum_r |\F{f}_s (r)|^{2s} \right)^{\frac{1}{2s}}
                    \cdot
$$
$$
                    \cdot
                    \left( \frac{1}{N} \sum_r |\F{g}_1 (r)|^{2t} \right)^{\frac{1}{2t}}
                        \dots
                            \left( \frac{1}{N} \sum_r |\F{g}_t (r)|^{2t} \right)^{\frac{1}{2t}}
                            =
$$
$$
                            =
                                (T_{s} (f_1))^{1/2s} \dots (T_{s} (f_{s}))^{1/2s}\,
                                (T_{t} (g_1))^{1/2t} \dots (T_{t} (g_{t}))^{1/2t} \,.
$$
This completes the proof.

\Cor
{\it
    Let $A,B$ be finite subsets of $G$.
    Then
    \begin{equation}\label{f:wrong_Minkovskii}
        T^{1/2k}_k (A\cup B)
            \le
                T^{1/2k}_k (A) + T^{1/2k}_k (B) \,.
    \end{equation}
}
\label{cor:wrong_Minkovskii}

We need in the notion of dissociativity in $\mathbf{F}^n_2$.

\Def
Let $R \subseteq \mathbf{F}^n_2$ be a set, $R=-R$ and $\{ 0 \} \in R$.
We say that a set $\L = \{ \l_1, \dots, \l_{|\L|} \} \subseteq \mathbf{F}^n_2$
belongs to the family $\mathbf{\L}_R (k)$
if the equality
   \begin{equation}\label{f:def_diss}
        \sum_{i=1}^{|\L|} \eps_i \l_i \in R \,,
   \end{equation}
where $\eps_i \in \{ -1,0,1 \}$ and $\sum_{i=1}^{|\L|} |\eps_i| \le k$
implies that all $\eps_i$ are equal to zero.
If $R=\{ 0\}$ then $\L$ belongs to the family $\mathbf{\L} (k)$.

\Pred
{\it
    Let $k$ be a positive integer, $k\ge 2$,
    and
    $\L \subseteq \mathbf{F}^n_2$ be an arbitrary set, belonging to the family $\mathbf{\L} (2k)$.
    Then for any integer $p$, $2\le p \le k$, we have
    \begin{equation}\label{est:T_p_L}
        T_p (\L) \le p^{p} |\L|^p \,.
    \end{equation}
}
\label{st:diss}
\Proof
Let $m = |\L|$.
Consider the equation
\begin{equation}\label{f:L_eq}
    \l_1 + \dots + \l_{2p} = 0, \quad \l_i \in \L, \quad i=1,\dots, 2p \,.
\end{equation}
Let us consider any partitions
$\mathcal{M} = \{ M_1, \dots, M_p \}$ of the segment $[2p]$
onto sets $M_j$, $|M_j| = 2$, $j=1,\dots, p$.
It is easy to see that  the number of such  partitions equals
$
    \frac{(2p)!}{2^p p!} \le \frac{(2p)^p}{2^p} = p^p
$.
Further, let us mark any set $M_j$ by an element $\l^{(j)}$ of the set $\L$.
Then the number of these {\it labelled} partitions does not exceed $p^p m^p$.
By assumption the set $\L$ belongs to the family $\mathbf{\L} (2k)$.
Hence if $(\l_1, \dots, \l_{2p})$ is an arbitrary solution of (\ref{f:L_eq})
then any $\l_i$, $i\in [2p]$
appears even number of times in this solution.
So a solution $(\l_1, \dots, \l_{2p})$ of (\ref{f:L_eq}) corresponds
a labelled partition $\mathcal{M}' = \{ (M_1, \l^{(1)}), \dots, (M_p, \l^{(p)}) \}$.
To see this let us construct
a labelled partition $\mathcal{M}' = \{ (M_1, \l^{(1)}), \dots, (M_p, \l^{(p)}) \}$
such that for any $M_j = \{ \a, \beta \}$, $j=1,\dots, p$, we have $\l_\a = \l_\beta = \l^{(j)}$.
Clearly,  if we have two different solutions of (\ref{f:L_eq})
then we get different labelled partitions.
Hence
the total number of solutions of (\ref{f:L_eq}) does not exceed $p^p m^p$.
This completes the proof.

\Note Rudin's Theorem (see \cite{Rudin_book,Rudin}) asserts
that for any functions $f:G \to \C$, $\supp \F{f} \subseteq \L$,
$\L$ is a dissociated set, we have $\| f \|_k \le C \sqrt{k} \| f \|_2$,
where $C>0$ is an absolute constant and $k\ge 2$.
In other words, for an arbitrary  $a_\l$ the following holds
\begin{equation}\label{f:note_1}
    \frac{1}{N} \sum_x \left| \sum_{\l\in \L} a_\l e(-\l\cdot x) \right|
        \le
            C^k k^{k/2} \left( \sum_{\l \in \L} |a_\l|^2 \right)^{k/2} \,.
\end{equation}
Certainly, inequality (\ref{f:note_1}) implies Proposition \ref{st:diss} :
to see this one can put $k=2p$ and $a_\l = 1$.
On the other hand, we can use a slightly modified arguments from Proposition \ref{st:diss}
to prove (\ref{f:note_1}).
Indeed, to obtain inequality (\ref{f:note_1}),
we need to
calculate
the number of solutions of (\ref{f:L_eq})
such that  any solutions has weight $a_{\l_1} \dots a_{\l_{2p}}$
By assumption the set $\L$ belongs to the family $\mathbf{\L} (2k)$.
Hence if $(\l_1, \dots, \l_{2p})$ is an arbitrary solution of (\ref{f:L_eq})
then any $\l_i$, $i\in [2p]$
appears even number of times in this solution.
It follows that if a partition
$\mathcal{M} = \{ M_1, \dots, M_p \}$ of the segment $[2p]$
onto the sets $M_j$, $|M_j| = 2$, $j=1,\dots, p$
is fixed
then we get weight $\left( \sum_{\l \in \L} |a_\l|^2 \right)^{p}$.
We know that the number of such partitions $\mathcal{M}$ does not exceed $p^p$.
Thus, we have proved (\ref{f:note_1}) in the case $k=2p$.
Using standard methods (see e.g.  \cite{Green_Chang2}, Lemma 19),
we
obtain inequality (\ref{f:note_1}) for all
$k\ge 2$.

Now we can prove an analog of Proposition \ref{st:diss}
for subsets of sums of dissociated sets
and obtain Theorem \ref{t:my_Bougrain_new}.

\Pred
{\it
    Let $k$, $d$ be positive integers, $k\ge 2$,
    and
    $\L \subseteq \mathbf{F}^n_2$ be an arbitrary set, $\L \in \mathbf{\L} (2dk)$
    such that $|\L| \ge 4d^2$.
    Let also $Q$ be a subset of $d\dot{\L}$.
    Then for all integer $p$, $2\le p \le k$, we have
    \begin{equation}\label{est:T_p_dL}
        T_p (Q) \le 2^{8dp} p^{dp} |Q|^p \,.
    \end{equation}
}
\label{st:dissd}
\Proof
We use induction.
If $d=1$, then the bound for $T_p (Q)$ was obtained in Proposition \ref{st:diss}.
Let $d\ge 2$, and let $m = |Q|$.
Put $c_d := 8d$, $d\ge 1$.

Let $a = [|\L|/2d]$.
By assumption $|\L| \ge 4d^2$.
Hence $|\L|/a \le 4d$.
Besides
$$
    \binom{|\L|-d}{a-1}^{-1} \binom{|\L|}{a}
        =
            \frac{|\L| (|\L|-1) \dots (|\L|-d+1)}{a (|\L|-a) (|\L|-a -1) \dots (|\L|-a-d+2)}
                \le
$$
\begin{equation}\label{f:zapas}
                \le
                    4d \cdot e^{\frac{1}{|\L|} (\sum_{i=1}^{d-1} i + 2 \sum_{i=0}^{d-2} (a+i) )}
                        \le
                            2^4 d \,.
\end{equation}
Let $E$ be any set.
By $E^c$ denote $\L \setminus E$.
Using dissociativity of $\L$
and the definition of the operation $\dotplus$, we get
$$
    Q(x)
        = d^{-1} \binom{|\L|-d}{a-1}^{-1}
            \sum_{\L_0 \subseteq \L,\, |\L_0| = a}
                \left( Q\bigcap (\L_0 + (d-1) \dot{\L}_0^c \right) (x) \,.
$$
Using H\"{o}lder's inequality, we obtain
\begin{equation}\label{f:16:42}
    T_p (Q)
        \le
            d^{-2p} \binom{|\L|-d}{a-1}^{-2p} \binom{|\L|}{a}^{2p-1}
                \sum_{\L_0 \subseteq \L,\, |\L_0| = a}
                    T_p ( Q\bigcap (\L_0 + (d-1) \dot{\L}_0^c ) \,.
\end{equation}
If we prove  for any $\L_0 \subseteq \L$
the following  inequality
$$
    T_p ( Q\bigcap (\L_0 + (d-1) \dot{\L}_0^c ) \le 2^{c_{d-1} p} p^{dp} |Q\bigcap (\L_0 + (d-1) \dot{\L}_0^c|^p \,,
$$
then using (\ref{f:16:42}) and (\ref{f:zapas}), we obtain
$$
    T_p (Q)
        \le
            d^{-2p} \binom{|\L|-d}{a-1}^{-2p} \binom{|\L|}{a}^{2p} 2^{c_{d-1} p} p^{pd} m^p
                \le 2^{(c_{d-1}+8)p} p^{pd} m^p
                = 2^{c_d p} p^{pd} m^p
$$
and Proposition \ref{st:dissd} is proved.

Let $\L_1 = \t{\L}$, $\L_2 = \L \setminus \t{\L}$,
and $Q' \subseteq \L_1 + (d-1) \dot{\L}_2$.
We have to prove that
$T_p (Q') \le 2^{c_{d-1}p} p^{pd} |Q'|^p$.
Let $\l$ be an element from $\L_1$.
Consider the sets
$$
    D_\l = D (\l) = \{~ \l' ~:~ \l \dotplus \l' \in Q',\, \l' \in (d-1)\dot{\L}_2 ~\} \,,
$$
$$
    Q_\l = Q (\l) = \{~ q\in Q' ~:~ q = \l \dotplus \l'_2 \dotplus \dots \dotplus \l'_d ,\, \quad \l'_i \in \L,\, i=2,\dots,d ~\} \,,
$$
Clearly,  $Q(\l) = D(\l) + \l$.
Let $s_1$ be a number of nonempty sets $D_\l$.
Let these sets are $D_{\l_1}, \dots, D_{\l_{s_1}}$.
We shall write $D_j$ instead of $D_{\l_j}$.
Let also $s_2 = |\L_2|$.
Obviously, $Q\subseteq \{ \l_1, \dots, \l_s \} + (d-1) \dot{\L}_2$

Consider the equation
\begin{equation}\label{f:1'}
    q_1 + \dots + q_p = q_{p+1} + \dots + q_{2p} \,,
\end{equation}
where $q_i \in Q'$, $i=1,\dots, 2p$.
Denote by $\sigma$ the number of solutions of (\ref{f:1'}).
Since $Q'\subseteq \L_1 + (d-1)\dot{\L}_2$, it follows that for all  $q\in Q'$, we have
$q=\l + \mu$, where $\l \in \L_1$, $\mu \in (d-1)\dot{\L}_2$.

Let $i_1,\dots, i_{2p} \in [s_1]$ be arbitrary numbers.
Denote by $\sigma_{\v{i}}$,
$\v{i} = (i_1,i_2, \dots, i_{2p})$
the set of solutions of equation (\ref{f:1'}) such that
for all $j\in [2p]$, we have the restriction
$q_{j} \in D ({\l_{i_j}})$, $\l_{i_j} \in \L_1$.
By assumption the sets $\L_1$ and $\L_2$ belong to the family $\mathbf{\L} (2dk)$.
Also $L_1$, $L_2$ have empty intersection.
It follows that if  $q_1, \dots, q_{2p}$ is a solution
of equation (\ref{f:1'}) such that this solution belongs the set  $\sigma_{\v{i}}$, then
any component of vector $\v{i}$ appears even number of times in the vector.
We have
\begin{equation}\label{f:2'}
    \sigma
        \le
            \sum_{\mathcal{M},\, \mathcal{M} = \{ M_1, \dots, M_r \},\,\, [2p] = M_1 \bigsqcup \dots \bigsqcup M_p}\,
                \sum_{\v{i} \in \mathcal{M}} | \sigma_{\v{i}} | \,.
\end{equation}
Summation in (\ref{f:2'}) is taken over families of sets
$\mathcal{M},\, \mathcal{M} = \{ M_1, \dots, M_p \},\, [2p] = M_1 \bigsqcup \dots \bigsqcup M_p$
such that for all $j=1,\dots, p$, we have $|M_j|=2$.
Let $M_j = \{ \a^{(j)}_1, \a^{(j)}_{2} \}$, $j=1,\dots,p$.
By definition $\v{i} \in \mathcal{M}$ if for all $j\in [p]$, we have
$i_{\a^{(j)}_1} = i_{\a^{(j)}_{2}}$.

Using Lemma \ref{l:conv} and induction, we get
$$
    |\sigma_{\v{i}}| \le 2^{c_{d-1}} p^{d(p-1)} \, \prod_{j=1}^{2p} |D(\l_{i_j})|^{1/2} \,.
$$
Hence
\begin{equation}\label{f:3'}
    \sigma
        \le
            2^{c_{d-1}} p^{d(p-1)} \sum_{\mathcal{M}}
                \sum_{\v{i} \in \mathcal{M}} \, \prod_{j=1}^{2p} |D(\l_{i_j})|^{1/2} \,.
\end{equation}
Let $m' = |Q'|$, and let $q$ be an arbitrary element of the set $Q'$.
By assumption $\L_1\bigcap \L_2 = \emptyset$ and $\L$ is a dissociated set,
so it is easy to see that the sets $Q(\l)$ are disjoint.
Hence
\begin{equation}\label{f:T}
    \sum_{\l \in \L_1} |D(\l)| = \sum_{\l \in \L_1} |Q(\l)| = m' \,.
\end{equation}
For any $\l\in \L_1$, we have
$|D_\l| \le m'$.
Let $x\ge 1$ be an arbitrary number.
Using formula (\ref{f:T}), we get
\begin{equation}\label{f:TT}
    \sum_{\l \in \L_1} |D(\l)|^x = \sum_{\l \in \L_1} |Q(\l)|^x \le (m')_2^{x-1} \sum_{\l \in \L_1} |Q(\l)| = (m')^x \,.
\end{equation}
The number of partitions $\mathcal{M}$ in inequality (\ref{f:3'}) does not exceed $p^p$.
Any component of vector $\v{i}$
appears even number of times in the vector.
Combining inequality (\ref{f:3'}) and bound (\ref{f:TT}), we obtain
$
    \sigma \le 2^{c_{d-1}} p^{dp} (m')^p
$.
This completes the proof.

In some sense the last proposition is best possible.

\Pred
{\it
    Let $k,d$ be positive integers, and a set $\L \subseteq \mathbf{F}_2^n$
    belongs to the family $\mathbf{\L} (2d)$.
    Let also $\L_1$ be an arbitrary subset of $\L$, and
    $Q = d \dot{\L}_1 \subseteq d \dot{\L}$.
    Then for all $k \le |\L_1| / (2d)$ and for any $2\le p \le k$, we have
    $
        T_p (Q) \ge 2^{-3pd} p^{pd} |Q|^p
    $.
}
\label{p:exact}
\\
\Proof
Consider the equation
\begin{equation}\label{JJ}
    q_1 + \dots + q_p = q_{p+1} + \dots + q_{2p} \,,
\end{equation}
where $q_i \in Q$, $i=1,\dots,2p$.
Let us prove that equation  (\ref{JJ}) has at least $2^{-3pd} p^{pd} |Q|^p$
solutions.
Since $q_i \in Q$, it follows that $q_i = \sum_{j=1}^d \l^{(i)}_{j}$, $i=1,\dots, 2p$.
Consider tuples $(q_1, \dots, q_p)$ such that {\it all} $\l^{(i)}_{j}$
for all $q_i$ are different.
Clearly, there are exactly
$
    \binom{|\L_1|}{pd} \frac{(pd)!}{(d!)^p}
$
of such tuples.
For any tuple $(q_1, \dots, q_p)$ there are at least $\frac{(pd)!}{(d!)^p}$
solutions of equation (\ref{JJ}).
Indeed we have $\frac{(pd)!}{(d!)^p}$ number of ways to partition the set
$\{ \l^{(1)}_1, \dots, \l^{(1)}_d, \dots, \l^{(p)}_1, \dots, \l^{(p)}_d \}
 = \{ \l_1, \dots, \l_{pd} \}$
onto $p$ sets $M_1, \dots, M_p$ of the same cardinality.
Put $q_i = \sum_{j\in M_i} \l_j$, $i=p+1,\dots,2p$, we get a tuple
$(q_{p+1}, \dots, q_{2p}) \in Q^p$ such that $q_1 + \dots + q_p = q_{p+1} + \dots + q_{2p}$.
By assumption $\L$ is a dissociated set, thus any collection of sets $M_1, \dots, M_p$
corresponds a tuple $(q_{p+1}, \dots, q_{2p})$.
Hence $T_p (Q) \ge \binom{|\L_1|}{pd} \frac{(pd)!}{(d!)^p} \cdot \frac{(pd)!}{(d!)^p}$.
Since $\L \in \mathbf{\L} (2d)$, it follows that $|Q| = \binom{|\L_1|}{d}$.
Using inequality $2kd \le |\L_1|$, we get
$$
    T_p (Q)
        \ge
            \binom{|\L_1|}{pd} \frac{(pd)!}{(d!)^p} \cdot \frac{(pd)!}{(d!)^p}
                \ge
                    2^{-pd} \frac{(pd)!}{(d!)^p} |Q|^p
                        \ge
                            2^{-3pd} p^{pd} |Q|^p \,.
$$
This completes the proof.

At the end of the section we show that
Theorem \ref{t:main} (actually Theorem \ref{t:main_G}, see Appendix)
and Proposition \ref{st:dissd}
imply Theorem \ref{t:Bougrain_new} in the case $G=\mathbf{F}_2^n$.

\Th
 \label{t:my_Bougrain_new}
 { \it
    Let $\d,\a$ be real numbers, $0<\a \le \d \le 1/4$,
    $d$ be a positive integer, $d\le \log (1/\d) / 4$,
    $A$ be an arbitrary subset of $\mathbf{F}_2^n$ of the cardinality $\d N$,
    and let $\r_\a$ as in (\ref{f:R_def}).
    Suppose that a set $\L \subseteq \mathbf{F}_2^n$
    belongs to the family $\mathbf{\L} (2 \log (1/\d))$.
    Then for all $1\le d \le \log (1/\d) / 4$, we have
    \begin{equation}\label{est:dL&R}
        |d\dot{\L} \bigcap \r_\a| \le \left( \frac{\d}{\a} \right)^2 \left( \frac{2^{12} \log (1/\d)}{d} \right)^d \,.
    \end{equation}
 }
\Proof
Let $k=[\ln (1/\d) / d] \ge 2$,
$Q = d\dot{\L} \bigcap \r_\a$
 and $m=|Q|$.
We need to prove that
$m\le ( \d/\a )^2 ( \frac{2^{12} \log (1/\d)}{d})^d$.
Using Theorem \ref{t:main_G},
we get
\begin{equation}\label{tMp}
    T_k (Q) \ge \frac{\d \a^{2k}}{\d^{2k}} m^{2k} \,.
\end{equation}
On the other hand, by Proposition \ref{st:dissd}, we obtain
$
   T_k (Q) \le 2^{8kd} k^{kd} m^k
$.
Combining the last inequality and  (\ref{tMp}), we get
$m\le ( \d/\a )^2 ( \frac{2^{12} \log (1/\d)}{d})^d$.
This concludes the proof.


So an upper bound for $|d\dot{\L} \bigcap \r_\a|$ was obtained in Theorem \ref{t:my_Bougrain_new}.
The next simple proposition gives us a lower estimate for the quantity.
It is turn out this lower bound is close to the upper one.

\Pred
\label{pr:my_Bougrain_new_low}
{\it
    Let $\d$ be a real number, $1/N \le \d \le 1/16$,
    and $\a = 2^{-12} \d /\sqrt{n}$, $n\ge 32$.
    Then there exist a set $A \subseteq \mathbf{F}_2^n$, $\d N \le |A| \le 8 \d N$
    and a dissociated set $\L \subseteq \r_\a (A)$ such that
    for all integers $d\ge 1$, we have
    \begin{equation}\label{}
        |d\dot{\L} \bigcap \r_\a| \ge 2^{-30} \left( \frac{\d}{\a} \right)^2 \left( \frac{\log (1/\d)}{16d} \right)^{d-1} \,.
    \end{equation}
}
\Proof
Let $\v{e}_1 = (1,0,\dots,0), \dots, \v{e}_n = (0,\dots,0,1)$ be the standard basis for $\mathbf{F}_2^n$.
Let also $k=[\log 1/(4\d)]$, and $H,H^\bot$ be subspaces spanned by
vectors $\v{e}_1, \dots, \v{e}_{n-k}$ and $\v{e}_{n-k+1}, \dots, \v{e}_{n}$, correspondingly.
Let $A\subseteq H$ be a set of $\v{x} = (x_1,\dots,x_n) \in H$
such that the number $x_i = 1$, $i=1,\dots, n-k$ is at least $(n-k)/2$.
Clearly, $|A| \ge 2^{n-k-2} \ge \d N$ and
$|A| \le |H| \le 2^{n-k} \le 8 \d N$.
Let $H'$ be a space spanned by vectors of the length $n-k$,
namely $(1,0,\dots,0), \dots, (0,\dots,0,1)$.
Let also $n' = n-k$, and let $A' \subseteq H'$ be the restriction of $A$ on $H'$.
Let us find Fourier coefficients of $A'$.
We have
\begin{equation}\label{f:(I)}
    \F{A}' (r)
        =
            \sum_x A' (x) (-1)^{<r, x>}
                =
                    |A' \bigcap H^{(0)}_r| - |A' \bigcap H^{(1)}_r|
                        =
                            2 |A' \bigcap H^{(0)}_r| - |A'| \,,
\end{equation}
where $H^{(0)}_r = \{ x\in H' ~:~ <r,x> =0 \}$ and $H^{(1)}_r = \{ x\in H' ~:~ <r,x> =1 \}$.
Let $l\ge 0$ be a positive integer.
Consider the sets
$$
    \mathcal{H}_l
        =
            \{ x = (x_1, \dots, x_{n'}) ~:~ \# x_i =1 \mbox{ equals } l \} \,.
$$
Let $r\in \mathcal{H}_1$.
Put $\binom{x}{y} = 0$ for $y>x$.
Using Stirling's formula and (\ref{f:(I)}), we get
$$
    |\F{A}' (r)|
        =
            \left| \sum_{s=\lceil n'/2\rceil}^{n'} \left( 2 \binom{n'-1}{s} - \binom{n'}{s} \right) \right|
                =
                    \sum_{s=\lceil n'/2\rceil}^{n'} \left( \frac{2s-n'}{n'} \right) \binom{n'}{s}
                        \ge
$$
\begin{equation}\label{f:(II)}
                        \ge
                            \sum_{s=\lceil n'/2 + \sqrt{n'}/2 \rceil}^{[n'/2 + \sqrt{n'}]}
                                \left( \frac{2s-n'}{n'} \right) \binom{n'}{s}
                                    \ge
                                        e^{-8} \frac{1}{2\sqrt{\pi}} \frac{2^{n'}}{\sqrt{n'}}
                                        \ge
                                            2^{-14}\, \frac{2^{n'}}{\sqrt{n'}}
                                                \ge
                                                    2^{-12}\, \frac{\d N}{\sqrt{n}} \,.
\end{equation}
It is easy to see that for any $r\in H'$ and for all $h^{\bot} \in H^{\bot}$,
we have $\F{A} (r+h^{\bot}) = \F{A}' (r)$.
Hence $\mathcal{H}_1 + H^{\bot} \subseteq \r_\a (A)$, $\a = 2^{-12} \d /\sqrt{n}$
and $|\r_\a (A)| \ge n' 2^k \ge n/(16\d) \ge 2^{-28} \cdot \d/\a^2$.
Thus we have a lower bound for the cardinality of $\r_\a (A)$, which
is close to an upper bound --- $\d/\a^2$.
Clearly, the set $A'$ is invariant under all permutations.
Using this fact one can prove (assuming some restrictions on parameters)
that the following holds $\r_\a (A) = (\{ 0 \} \sqcup \mathcal{H}_1) + H^{\bot}$.
We do not need in the fact.

Let $\L = \{ \v{e}_1, \dots, \v{e}_{n} \} \subseteq \r_\a (A)$,
and $\L^* = \{ \v{e}_{n-k+1}, \dots, \v{e}_{n} \}$.
Clearly, \\
$\bigsqcup_{h_1 \in \mathcal{H}_1} (h_1 + (d-1) \dot{\L}^*) \subseteq \r_\a (A) \bigcap d \dot{\L}$.
Hence
$$
    |\r_\a (A) \bigcap d \dot{\L}|
        \ge
            n' \binom{k}{d-1}
                \ge
                    \frac{n}{4} \cdot \frac{k^{d-1}}{d^{d-1} e^{d-1}}
                        \ge
                            2^{-30} \left( \frac{\d}{\a} \right)^2 \left( \frac{\log (1/\d)}{16d} \right)^{d-1} \,.
$$
This completes the proof.

\Note Certainly, we can change the value of the parameter $\a$
in Proposition \ref{pr:my_Bougrain_new_low}.
For example one can consider sets $\mathcal{H}_2$
instead of $\mathcal{H}_1$
and choose the parameter $\a$ smaller than $2^{-12} \d /\sqrt{n}$.

\refstepcounter{section}
\label{connected_sets}

\begin{center}
{\large
{\bf \arabic{section}. On connected subsets of $d\dot{\L}$.}}
\end{center}

Let $G$ be an Abelian group, and $A \subseteq G$ be an  arbitrary finite set.
In paper \cite{Shkr_small}, so--called "connected"\, sets $A$
were studied  (see also article \cite{Ruzsa_Elekes}).
Let us give a definition from \cite{Shkr_small}.

\Def
\label{b_connectedness+}
    Let $k\ge 2$ be a positive integer,
    and $\beta_1, \beta_2 \in [0,1]$ be real numbers, $\beta_1 \le \beta_2$.
    Nonempty set $A\subseteq G$ is called
    {\it $(\beta_1, \beta_2)$--connected of degree $k$} if
    there exists an absolute constant $C\in (0,1]$
    such that for any $B\subseteq A$,
    $\beta_1 |A| \le |B| \le \beta_2 |A|$ we have
    \begin{equation}\label{ineq:b_connectedness+}
        T_k (B) \ge C^{2k} \left( \frac{|B|}{|A|} \right)^{2k} T_k (A) \,.
    \end{equation}

By $\zeta_k (A)$ denote the quantity
$
    \zeta_k (A) := \frac{\log T_k (A)}{\log |A|}
$.
In paper \cite{Shkr_small} (see also \cite{Sanders_Sh})
the following result was obtained.

\Th
{\it
    Let $\beta_1, \beta_2 \in (0,1)$ be real numbers, $\beta_1 \le \beta_2$.
    Then there exists a set $A' \subseteq A$ such that \\
    $1)~$ $A'$ is $(\beta_1, \beta_2)$--connected of degree $k$ set
            such that (\ref{ineq:b_connectedness+}) holds for any $C \le 1/32$.\\
    $2)~$ $|A'| \ge m \cdot 2^{\frac{\log ( (2k-1)/ \zeta)}{\log (1+\kappa)} \log (1-\beta_2)}$,
                    where $\kappa = \frac{\log ((1-\beta_1)^{-1})}{\log m} (1-16 C)$.\\
    $3)~$ $\zeta_k (A') \ge \zeta_k (A)$.
}
\label{t:connected_old}

In the section we prove an analog of Theorem \ref{t:connected_old}
for subsets of dissociated sets.

Let $\L \subseteq \mathbf{F}^n_2$ be an arbitrary set from the family $\mathbf{\L} (2dk)$,
and $A\subseteq d\dot{\L}$.
    Denote by $D_k (A)$ the quantity
    \begin{equation}\label{}
        D_k (A) = \log \left( \frac{T_k (A)}{k^k |A|^k} \right) \,.
    \end{equation}
    In other words $T_k (A) = 2^{D_k (A)} k^k |A|^k$.
    Since for all sets $A$ with sufficiently large cardinality,
    we have $T_k (A) \ge \binom{|A|}{k} (k!)^2 \ge e^{-2k} k^k |A|^k$, it follows that
    the quantity $D_k (A)$ is at least
    $-2 k \log e$.
    On the other hand, by Proposition \ref{st:dissd}, we get
    $D_k (A) \le 8 d \log d + k (d-1) \log k$.

\Th
{\it
    Let $K>0$ be a real number, $k,d$ be positive integers, $k,d\ge 2$.
    Let $\L \subseteq \mathbf{F}^n_2$ be an arbitrary set, $\L \in \mathbf{\L} (2dk)$,
    and $Q$ be a subset of $d\dot{\L}$ such that
    $T_k (Q) \ge \frac{k^{dk} |Q|^k}{K^k}$.
    Let also $\beta_1, \beta_2 \in (0,1)$ be real numbers, $\beta_1 \le \beta_2$.
    Then there is a set $Q' \subseteq Q$ such that \\
    $1)~$ $Q'$ is a $(\beta_1, \beta_2)$--connected of degree $k$
               such that (\ref{ineq:b_connectedness+}) holds for any $C\le 1/8$.\\
    $2)~$ $|Q'| \ge |Q| \cdot 2^{\frac{8d \log d + k (d-1) \log k - D_k (Q)}{k \log (1+\beta_1 (1-4C)) } \log (1-\beta_2)}$.\\
    $3)~$ $T_k (Q') \ge \frac{k^{dk} |Q'|^k}{K^k}$.
}
\label{t:connected}
\\
\Proof
    Let $m=|Q|$, and $C\le 1/8$ be a real number.
    The proof of Theorem \ref{t:connected} is a sort of algorithm.
    If $Q$ is $(\beta_1, \beta_2)$--connected of degree $k$
    and (\ref{ineq:b_connectedness+}) is true with the constant $C$
    then there is nothing to prove.
    Suppose that $Q$ is not
    $(\beta_1, \beta_2)$--connected of degree $k$ set (with the constant $C$).
    Then there exists a set $B\subseteq Q$, $\beta_1 |Q| \le |B| \le \beta_2 |Q|$
    such that (\ref{ineq:b_connectedness+}) does not hold.
    Note that $|Q|>2$.
    Let $\ov{B} = Q\setminus B$ and $c_B = |B|/ |Q|$.
    We have $\beta_1 \le c_B \le \beta_2$.
    Using Corollary \ref{cor:wrong_Minkovskii}, we get
    \begin{equation}\label{I:f}
        T_k (\ov{B}) > T_k (Q) (1-C c_B)^{2k} \,.
    \end{equation}
    Let $b = |B|$ and $\ov{b} = |\ov{B} | = m - b$, $D=D_k (Q)$, $\ov{D} = D_k (\ov{B})$.
    By inequality (\ref{I:f}), we obtain
$$
    \ov{D}
        >
            D + k\log m - k \log \ov{b} + 2k \log (1-Cc_B)
                =
                    D + k \left(\log (\frac{m}{m-b} (1-Cc_B)^2) \right)
                        \ge
$$
\begin{equation}\label{f:t_1.3}
                        \ge
                            D + k\log ((1+c_B)(1-2Cc_B))
                            \ge
                                D + k \log (1+\beta_1 (1-4C)) \,.
\end{equation}
    Besides, by the definition of  $(\beta_1, \beta_2)$--connectedness of degree $k$,
    we have
    \begin{equation}\label{f:t_1.4}
        |\ov{B}| \ge (1-\beta_2) m = (1-\beta_2) |Q| \,.
    \end{equation}
    Thus if the set $Q$ is not $(\beta_1, \beta_2)$--connected of degree $k$
    then there is a set $\ov{B} \subseteq Q$ such that
    (\ref{f:t_1.3}), (\ref{f:t_1.4}) hold.
    Put $Q_1 = \ov{B}$ and apply the arguments above to $Q_1$.
    And so on.
    We get the sets $Q_0 = Q, Q_1, Q_2, \dots, Q_s$.
    Clearly, for any  $Q_i$, we have $D_k (Q_i) \le 8 d \log d + k (d-1) \log k$.
    Using this and (\ref{f:t_1.3}), we obtain that the total number of steps
    of our algorithm
    does not exceed
    $
    \frac{8d \log d + k (d-1) \log k - D_k (Q)}{k \log (1+\beta_1 (1-4C)) }
    $.
    At the last step of the algorithm, we find the set $Q'=Q_s \subseteq Q$
    such that $Q'$ is $(\beta_1, \beta_2)$--connected of degree $k$
    and such that $D_k (Q') \ge D_k (Q)$.
    Thus  $Q'$ has the properties $1)$ and $3)$
    of the  Theorem.
    Let us prove $2)$.
    Using (\ref{f:t_1.4}), we obtain
    $$
        |Q'|
            \ge
                (1-\beta_2)^s m
                    \ge
                            m \cdot 2^{\frac{8d \log d + k (d-1) \log k - D_k (Q)}{k \log (1+\beta_1 (1-4C)) } \log (1-\beta_2)} \,.
    $$
    This concludes the proof.

    We shall use Theorem \ref{t:connected} in the next section.

\refstepcounter{section}
\label{inverse}

\begin{center}
{\large
{\bf \arabic{section}. On large subsets of sum of two dissociated sets.}}
\end{center}

Let $H = (h_{ij})$ be a matrix of the size $x\m y$, $x\le y$.
By $\per H$ denote the permanent of matrix $H$.
Recall that
\begin{equation}\label{f:permanent}
    \per H = \sum_{\sigma} h_{1\sigma(1)} \dots h_{x\sigma(x)} \,,
\end{equation}
where the summation in (\ref{f:permanent})
is taken over all injective maps $\sigma : [x] \to [y]$.
We need in a well--known Frobenius--K\"{o}nig's Theorem on nonnegative matrices (see \cite{Mink}).

\Th
{\it
    Let $p$ and $r$ be positive integers, $r\le p$, and
    let $H$ be a nonnegative matrix of size $p\m r$.
    Then the permanent of matrix $H$ equals zero iff
    $H$ contains a zero matrix of size $p-s+1\m s$.
}
\label{t:Frobenius-Konig}

Using Theorem \ref{t:Frobenius-Konig}, we prove
a simple lemma.

\Lemma
{\it
    Let $p$ and $r$ be positive integers,
    and let $H = (h_{ij})$ be a nonnegative matrix of size $p\m r$.
    Let also \\
    $1)~$ For all $i\in [p]$, we have $\sum_{j=1}^r h_{ij} \ge 2$.\\
    $2)~$ For all $j\in [r]$, we have $\sum_{i=1}^p h_{ij} \ge 1$, and,
    finally,\\
    $3)~$ $\sum_{i=1}^p \sum_{j=1}^r h_{ij} = 2p$.\\
    Deleting from $H$ all columns such that $\sum_{i=1}^p h_{ij} = 1$,
    we get matrix  $H_0$.
    Then the permanent of
    $H_0$
    does not equal zero.
}
\label{l:per>0}
\\
\Proof
Let the number of $j$ such that  $\sum_{i=1}^p h_{ij} = 1$ equals $e$.
Without loss of generality we can suppose that the matrix $H_0$
was obtain from $H$ by deletion of the last $e$ columns.
Let $H_0 = (h^{0}_{ij})$, $i=1,\dots,p$, $j=1,\dots, r-e = r_0$.
Applying condition $3)$ of the Lemma, we get
$\sum_{i=1}^p \sum_{j=1}^{r_0} h^{0}_{ij} = 2p-e$.
Using condition $2)$, we obtain $r_0 \le p$.
Suppose that our Lemma is false.
If the permanent of
$H_0$ equals zero then by Theorem \ref{t:Frobenius-Konig}
the matrix contains a submatrix of the size $s\m t$, $s+t=p+1$.
Using permutations of rows and columns, we can suppose that $H_0$
is
\begin{displaymath}
    H_0 =
        \left( \begin{array}{cc}
                    \mathbf{X} & \mathbf{Z} \\
                    \mathbf{0} & \mathbf{Y}
        \end{array} \right) \,,
\end{displaymath}
where zero matrix $\mathbf{0}$ has the size $s\m t$, $s+t=p+1$.
Denote by $s_1$ the number of $i\in \{ p-s+1, \dots, p \}$ such that
$\sum_{j=1}^{r_0} h^{0}_{ij} = 1$,
and by $s_2$ the number of $i\in \{ p-s+1, \dots, p \}$ such that
$\sum_{j=1}^{r_0} h^{0}_{ij} \ge 2$.
Clearly,
$s_1 \le e$.
By $2)$, we obtain $s=s_1 + s_2$.
Using condition $1)$ of the lemma, we get
$$
    2p-e = \sum_{i=1}^p \sum_{j=1}^{r_0} h^{0}_{ij}
            \ge
                \sum_{j=1}^t \sum_{i=1}^p h^{0}_{ij} + \sum_{i=p-s+1}^p \sum_{j=1}^{r_0} h^{0}_{ij}
                    \ge 2t + s_1 + 2s_2 = 2t + 2s - s_1 = 2p + 2 - s_1 \,.
$$
The last inequality implies $s_1 \ge e+2$
with contradiction.
This completes the proof.

Let $p$ be a positive integer,
$\L \subseteq \mathbf{F}^n_2$ be an arbitrary set, $\L \in \mathbf{\L} (2p)$,
and $\mathcal{E} = \{ E_1,\dots, E_{2p} \}$ be a tuple of subsets of $\L$.
In the proof of Proposition
\ref{st:dissd}
we estimated the number of solutions of the equation
\begin{equation}\label{f:explain}
    \l_1+\dots + \l_{2p} = 0,\quad \mbox{ where } \quad \l_i \in E_i,\, \quad i=1,\dots,2p \,.
\end{equation}
To calculate the number of such solutions,
we used Lemma \ref{l:T_k(B,A)}  --- a simple corollary
of H\"{o}lder inequality.
In the proof of the main result of this section --- Theorem \ref{t:inverse2},
we need in a more delicate result on the number of
solutions of equation (\ref{f:explain}).

\Lemma
{\it
    Let $p$ a positive integer,
    $\L  \subseteq \mathbf{F}^n_2$ be an arbitrary set, $\L \in \mathbf{\L} (2p)$,
    and $\mathcal{E} = \{ E_1,\dots, E_{2p} \}$ be a tuple of subsets of $\L$.
    Suppose that we have a partition of the segment $[2p]$ onto $r$ classes $\mathcal{C}_1,\dots,\mathcal{C}_r$.
    Let $S^* \subseteq [2p]$ be an arbitrary set, and $\bar{S}^* = [2p]\setminus S$.
    Let also $M(S^*) = (m_{ij})$ be a matrix of size $p\m p$,
    $m_{ij} = |E_i \bigcap E_j|$, $i\in S$, $j\in \bar{S}^*$.
    Then number of solutions of the equation
    \begin{equation}\label{f:E_i}
        \l_1 + \dots + \l_{2p} = 0,\quad \mbox{ where } \quad \l_i \in E_i,\, \quad i=1,\dots,2p
    \end{equation}
    does not exceed
    \begin{equation}\label{f:s(E)}
                \sum_{S^* \subseteq [2p],\, |S^*| = p}\,\, \per M(S^*) \,,
    \end{equation}
    and the summation in formula (\ref{f:s(E)}) is taken over
    all sets $S^*$ such that
    $S^*$ contains an element
    from any class
    $\mathcal{C}_i$ such that $|\mathcal{C}_i| \ge 2$.
}
\label{l:sophisticated}
\\
\Proof
Denote by $Z$ the number of solutions of equation (\ref{f:E_i}).
By assumption the set $\L$ belongs to the family $\mathbf{\L} (2p)$.
Hence
if $(\l_1, \dots, \l_{2p})$ is an arbitrary solution of (\ref{f:E_i})
then any $\l_i$, $i\in [2p]$ appears even number of times in this solution.
Thus (see proof of Proposition \ref{st:diss}), we get
\begin{equation}\label{f:rough}
    Z
        \le
            \sum_{\mathcal{K},\, \mathcal{K} = \{ K_1, \dots, K_p \},\,\, [2p] = K_1 \bigsqcup \dots \bigsqcup K_p}\,
                \prod_{j=1}^p \left| \bigcap_{\a \in K_j} E_{\a} \right| = Z_1 \,.
\end{equation}
The summation in (\ref{f:rough}) is taken over families of sets
$\mathcal{K},\, \mathcal{K} = \{ K_1, \dots, K_p \},\, [2p] = K_1 \bigsqcup \dots \bigsqcup K_p$
such that for any $j\in [p]$, we have $|K_j| = 2$.
Let us prove that
\begin{equation}\label{f:s(E)_1}
    Z_1 \le \sum_{S^* \subseteq [2p],\,|S^*| = p}\,\, \per M(S^*) \,,
\end{equation}
and the summation in formula (\ref{f:s(E)_1}) is taken over all sets $S^*$ such that
    $S^*$ contains an element
    from any class
    $\mathcal{C}_i$ such that $|\mathcal{C}_i| \ge 2$.
Clearly,
\begin{equation}\label{f:s(E)_2}
    \sum_{\mathcal{K},\, \mathcal{K} = \{ K_1, \dots, K_p \},\,\, [2p] = K_1 \bigsqcup \dots \bigsqcup K_p}\,
                \prod_{j=1}^p \left| \bigcap_{\a \in K_j} E_{\a} \right|
                    \le
                        \sum_{S^* \subseteq [2p],\,|S^*| = p}\,\, \per M(S^*)
\end{equation}
Indeed if $x$ is a summand from the left hand side of (\ref{f:s(E)_2})
which corresponds some partition $\mathcal{K}$
then $x$ is present in the right hand side too.
To see this let $S^*$ be the set of all first elements of $K_j$, $j=1,\dots,p$.
Let $\a$ is any of such numbers, $\a \in K_j$.
Then  there is quantity $|E_\a \bigcap E_\beta|$ with $\beta \in K_j$ in the right hand side of (\ref{f:s(E)_2}).
Taking  a product of such quantitaes, we get $x$.
Further if $y$ is an arbitrary summand from the right hand side of (\ref{f:s(E)_2})
then it is easy to form a partition $\mathcal{K}$ which corresponds to the $y$.

If $x$ is a summand from the left hand side of (\ref{f:s(E)_1})
which corresponds some partition $\mathcal{K}$
and we shall find a set $S^*$ such that
for all $j\in [p]$, we have $|K_j \bigcap S^*| = 1$
and such that $S^*$ contains an element from any class
$\mathcal{C}_i$ with restriction $|\mathcal{C}_i| \ge 2$
then we shall prove (\ref{f:s(E)_1}).
Let $H = (h_{\gamma \d})$ be a nonnegative matrix $p\m r$
such that any element $h_{\gamma \d}$ of $H$ equals $|K_\gamma \bigcap \mathcal{C}_\d|$.
Clearly, for all $\gamma \in [p]$, we have $\sum_\d h_{\gamma \d} = |K_\gamma| = 2$
and $\sum_{\gamma,\d} h_{\gamma \d} = \sum_{\gamma} |K_\gamma| = 2p$.
Since the sets $\mathcal{C}_1,\dots,\mathcal{C}_r$ form a partition of the segment $[2p]$,
it follows that for all $\d \in [r]$ the following holds
$\sum_\gamma h_{\gamma \d} = |\mathcal{C}_\d| \ge 1$.
Using Lemma \ref{l:per>0}, we obtain that the permanent of the matrix $H_0$
does not  equal zero.
Hence
$H_0$ contains a diagonal of nonzero elements.
Let the size of $H_0$ be $p\m r_0$.
Without loss of generality we can assume that the matrix $H_0$
is formed by first $r_0$ columnes of $H$.
Then nonzero diagonal $H_0$ is
$(\gamma_1, 1), \dots, (\gamma_{r_0}, r_0)$
and
for any $i\in [r]$ there is a number $\a_i \in K_{\gamma_i}$
such that $\a_i \in \mathcal{C}_i$ and $|\mathcal{C}_i| \ge 2$.
Let us add elements $\a_1, \dots, \a_r$ into the set $S^*$.
Besides let us add the first elements of all
$K_{\gamma}$, $\gamma \neq \gamma_1, \dots, \gamma_r$ in $S^*$.
It is easy to see that
$S^*$ contains an element
from any class $\mathcal{C}_i$ such that $|\mathcal{C}_i| \ge 2$.
This completes the proof.

\Note Lemma \ref{l:sophisticated} gives us an upper bound for $T_p (E_1,\dots,E_{2p})$.
This estimate implies (up to constants) the
bound $p^p \prod_{\a=1}^{2p} |E_\a|^{1/2}$
for $T_p (E_1,\dots,E_{2p})$,
which can be derived from Lemma \ref{l:conv}.
Indeed for any sets $A$ and $B$, we get
\begin{equation}\label{f:wrong_mixing}
    |A\bigcap B| \le \min\{ |A|, |B| \} \le |A|^{1/2} |B|^{1/2} \,.
\end{equation}
So any summand in $\per M(S^*)$ does not exceed $\prod_{\a=1}^{2p} |E_\a|^{1/2}$.
We have exactly $p!$ of such summands.
Thus by Lemma \ref{l:sophisticated}, we obtain  that
$T_p (E_1,\dots,E_{2p}) \le 2^{2p} p! \prod_{\a=1}^{2p} |E_\a|^{1/2}$.

\Lemma
{\it
    Let $\d_0 > 0$ be a real number, $r,p$ be positive integers,
    $p\ge 2\d_0 + 3$, $r\ge p - \d_0$.
    Let $t_1,\dots,t_r$ be a sequence of natural numbers such that
    $t_j \ge 2$, $j=1,\dots,r$ and $\sum_{j=1}^r t_j = 2p$.
    Let also $T = \max_{j\in [r]} t_j$, and
    $\a_j = | \{ j \in [r] ~:~ t_j \ge T-i \} |$, $i=0,1,\dots, T-2$.
    Let $z$ be a nonnegative number such that
    $\sum_{i=0}^{z-1} \a_i \le p < \sum_{i=0}^z \a_i$, and let $q_z = p - \sum_{i=0}^{z-1} \a_i$.
    Then the quantitaty
    $$
        \pi (t_1, \dots, t_r) := T^{\a_0} (T-1)^{\a_1} \dots (T-(z-1))^{\a_{z-1}} (T-z)^{q_z}
    $$
    does not exceed $ 2^{3p} \max\{\, \d_0^{4\d_0}, 1 \,\}$.
}
\label{l:columns}
\\
\Proof
Suppose that $\d_0 \ge 1$.
It is easy to see that the sequence
$\a_0, \a_1, \dots, \a_{T-2}$ is nondecreasing
and $\sum_{i=0}^{T-2} \a_i = \sum_{j=1}^r t_j = 2p$.
Using the condition $\sum_{j=1}^r t_j = 2p$ one more and inequalities $r\ge p-\d_0$,
$t_j \ge 2$, $j\in [r]$, we get
$T + 2(r-1) \le 2p$ and $T \le 2 \d_0 + 2 \le 4 \d_0$.
Suppose that $\a_0 = \dots = \a_{z-1} = q_z = 1$.
Then $p = \sum_{i=0}^{z-1} \a_i + q = z+1$.
On the other hand, there are exactly $T-1$ numbers $\a_i$.
Hence $z$ does not exceed $T-1$ and we get inequality $p\le T \le 2\d_0 + 2$
with contradiction.
Thus either $\a_{z-1} \ge 2$ or $q_z \ge 2$.

Let $\pi^*$ be the maxiamal value of the function $\pi (t_1, \dots, t_r)$
such that all $t_i$ satisfy
\begin{equation}\label{cond:admissible}
    t_1 + \dots + t_r = 2p, \quad t_j \ge 2, \quad r\ge p - \d_0 \,.
\end{equation}
If (\ref{cond:admissible}) holds for a tuple $t_1, \dots, t_r$
then we shall say that this tuple is {\it admissible}.
Let $\pi^* = \pi (t^0_1, \dots, t^0_r)$.
Without loss of generality we can assume that
$t^0_1 \ge t^0_2 \ge \dots \ge t^0_r$.
We have that either $\a_{z-1} \ge 2$ or $q_z \ge 2$.
Suppose that $t^0_2 \ge 3$.
Then we can consider an admissible tuple
$\t{t}_1 = t^0_1 + 1$, $\t{t}_2 = t^0_1 - 1$,
$\t{t}_3 = t^0_3, \dots, \t{t}_r = t^0_r$.
Clearly, $\pi^* = \pi (t^0_1, \dots, t^0_r) < \pi (\t{t}_1, \dots, \t{t}_r)$.
Whence $t^0_2 = 2$ and $\pi^* \le T^T 2^p \le 2^{3p} \d_0^{4\d_0}$.

Now suppose that $\d_0 < 1$.
In the case we have $T\le 4$.
Using a trivial estimate
$\pi (t_1, \dots, t_r) \le T^p \le 2^{2p}$
we get the required result.
This completes the proof.

Let $k$ be positive integer, $k\ge 2$, and
$\L_1,\L_2 \subseteq \mathbf{F}^n_2$ be arbitrary {\it disjoint} sets
such that $\L_1 \bigsqcup \L_2$ belongs to the family $\mathbf{\L} (4k)$.
Let also $Q$ be a subset of $\L_1 \dotplus \L_2 = \L_1 + \L_2$.
Define the sets $D (\l) = D_\l$ and $Q (\l) = Q_\l$, $\l\in \L_1$
(see proof of Proposition \ref{st:dissd}).
Let $\l \in \L_1$ and
$$
    D (\l) = \{~ \l' ~:~ \l + \l' \in Q,\, \l' \in \L_2 ~\} \,,
$$
$$
    Q (\l) = \{~ q\in Q ~:~ q = \l + \l' ,\, \l' \in \L_2 ~\} \,.
$$
Clearly, $Q(\l) = D(\l) + \l$.
Let $s_1$ be a number of nonempty sets $D_\l$.
Let these sets are $D_{\l_1}, \dots, D_{\l_{s_1}}$.
We shall write $D_j$ instead of $D_{\l_j}$.
Let also $s_2 = |\L_2|$.
Obviously, $Q\subseteq \{ \l_1, \dots, \l_s \} + \L_2$.

\Pred
{\it
    Let $M>0$ be a real number, $p$ be a positive integer, $p\ge 5$,
    and
    $\L_1,\L_2 \subseteq \mathbf{F}^n_2$ be arbitrary disjoint sets from the family $\mathbf{\L} (4p)$.
    Let also $Q$ be a subset of $\L_1 + \L_2$,
    $|Q| \ge \max\{ 2s_2 p, 2^8 s_2 p M^8 \}$,
    $\d_0 = \max\{ (p \log (2 e M) ) / \log (|Q|/(s_2 p)) , 1 \}$,
 and $X = \max\{\, \d_0^{4\d_0}, 1 \,\}$.
    Then
    \begin{equation}\label{e:T_p_intesect}
        T_p (Q)
            \le
                2^{5p} X p^{3p} s_2^p
                \cdot
                \sum_{r=p-\lceil \d_0 \rceil}^{p}
                \left( \frac{1}{ps_2} \right)^r
                    \cdot
                        \left( \sum_{S\subseteq [s_1],\, |S| = r}\, \prod_{\a\in S}
                            \left( \sum_{\beta \in S} |D_\a \bigcap D_\beta| \right) \right)
                                + \frac{p^{2p} |Q|^p}{2M^p} \,.
    \end{equation}
}
\label{st:inverse2}
\Proof
Let $m = |Q|$.
Consider the equation
\begin{equation}\label{f:I}
    q_1 + \dots + q_{2p} = 0 \,,
\end{equation}
where $q_i \in Q$, $i=1,\dots, 2p$.
Denote by $\sigma$ the number of solutions of equation (\ref{f:I}).
Since $Q\subseteq \L_1 + \L_2$, it follows that for all $q\in Q$, we have
$q=\l_1 + \l_2$, where $\l_1 \in \L_1$, $\l_2 \in \L_2$.

Let $i_1,\dots, i_{2p} \in [s_1]$ be arbitrary numbers.
By
$\sigma_{\v{i}}$,
$\v{i} = (i_1,i_2, \dots, i_{2p})$
denote the set of
solutions of equation (\ref{f:I}) such that
for all $j\in [2p]$ we have the restriction
$q_{j} \in D ({\l_{i_j}})$, $\l_{i_j} \in \L_1$.
By assumption the set $\L_1 \bigsqcup \L_2$ belongs to the family $\mathbf{\L} (4k)$
and $\L_1 \bigcap \L_2 = \emptyset$.
Hence if $(q_1, \dots, q_{2p}) \in \sigma_{\v{i}}$ is an arbitrary solution
of (\ref{f:I}) then any component of vector $\v{i}$ appears even number
of times in this vector.
We have
\begin{equation}\label{f:II}
    \sigma
        \le
            \sum_{\mathcal{N},\, \mathcal{N} = \{ N_1, \dots, N_r \},\,\, [2p] = N_1 \bigsqcup \dots \bigsqcup N_r}\,
                \sum_{\v{i} \in \mathcal{N}} | \sigma_{\v{i}} | \,.
\end{equation}
The summation in the right hand side of (\ref{f:II}) is taken over families of
sets $\mathcal{N},\, \mathcal{N} = \{ N_1, \dots, N_r \},\, [2p] = N_1 \bigsqcup \dots \bigsqcup N_r$
such that for all $j=1,\dots, r$ the cardinality of $N_j$ is an even number and $|N_j| \ge 2$.
Let $N_j = \{ \a^{(j)}_1, \dots, \a^{(j)}_{|N_j|} \}$, $j=1,\dots,r$.
By definition $\v{i} \in \mathcal{N}$ if for all $j\in [r]$ the following holds
$i_{\a^{(j)}_1} = \dots = i_{\a^{(j)}_{|N_j|}}$
and for any different sets $N_{j_1}$, $N_{j_2}$ from the partition $\mathcal{N}$,
we have $i_{\a} \neq i_{\beta}$, where $\a$ is an arbitrary element from $N_{j_1}$,
and $\beta$ is an element from $N_{j_2}$.

By $r=r(\mathcal{N})$ denotes the number of the sets $N_j$ in the partition $\mathcal{N}$.
We have
\begin{equation}\label{f:III}
    \sigma
        =
            \sum_{r=0}^{p-\lceil \d_0\rceil}
                \sum_{\mathcal{N},\, r(\mathcal{N}) = r}\,
                    \sum_{\v{i} \in \mathcal{N}} | \sigma_{\v{i}} |
                        +
            \sum_{r=p-\lceil \d_0 \rceil+1}^p
                \sum_{\mathcal{N},\, r(\mathcal{N}) = r}\,
                    \sum_{\v{i} \in \mathcal{N}} | \sigma_{\v{i}} | = \sigma_1 + \sigma_2 \,.
\end{equation}
Let us estimate the sum $\sigma_1$.
Let $q$ be an arbitrary element of the set $Q$.
Using the condition $\L_1\bigcap \L_2 = \emptyset$ and dissociativity of $\L$,
it is easy to see that the sets $Q(\l)$ are disjoint.
Hence
\begin{equation}\label{f:t}
    \sum_{\l \in \L_1} |D(\l)| = \sum_{\l \in \L_1} |Q(\l)| = m \,.
\end{equation}
For any $\l\in \L_1$, we have $|D_\l| \le s_2$.
Let $x\ge 1$ be an arbitrary number.
Using formula (\ref{f:t}), we get
\begin{equation}\label{f:tt}
    \sum_{\l \in \L_1} |D(\l)|^x = \sum_{\l \in \L_1} |Q(\l)|^x \le s_2^{x-1} \sum_{\l \in \L_1} |Q(\l)| = s_2^{x-1} m \,.
\end{equation}
Let $S_{\v{i}} = \{ i_j \}_{j\in [2p]}$.
Applying Lemma \ref{l:conv} and inequality (\ref{f:tt}), we get
$$
    \sigma_1
        \le
            \sum_{r=0}^{p-\lceil \d_0 \rceil}
                \sum_{\mathcal{N},\, r(\mathcal{N}) = r}\,
                    \sum_{\v{i} \in \mathcal{N}} \prod_{\a\in [2p]} |D_{i_\a}|^{1/2}
                        \le
                            p^p \sum_{r=0}^{p-\lceil \d_0 \rceil}
                                 \sum_{\mathcal{N},\, r(\mathcal{N}) = r}\,
                                    s_2^{p-r} \sum_{\v{i} \in \mathcal{N}} \prod_{\a \in S_{\v{i}}} |D_\a| \,.
$$
Note that if the lengths of the sets $N_j$ are fixed then the set $S_{\v{i}}$
does not change after any permutation of these sets.
Let $t_j = |N_j|$, $j=1,\dots, r$.
Using inequality $m \ge 2s_2 p$, the definition of the quantity $\d_0$
and identity (\ref{f:t}), we obtain
$$
    \sigma_1
        \le
            p^p \sum_{r=0}^{p-\lceil \d_0 \rceil}
                \sum_{t_1+\dots+t_r = 2p} \frac{(2p)!}{t_1! \dots t_r!}\, \frac{1}{r!}\, s_2^{p-r}\, m^r
                    \le
                        e^p p^p s_2^p \sum_{r=0}^{p-\lceil \d_0 \rceil} \left( \frac{m}{s_2} \right)^r r^{2p-r}
                            \le
$$
\begin{equation}\label{f:final1}
                            \le
                                2 e^p p^{3p} s^p_2 \left( \frac{m}{p s_2} \right)^{p-\d_0}
                                    =
                                        2 e^p p^{2p} m^p \left( \frac{s_2 p}{m} \right)^{\d_0}
                                            \le \frac{p^{2p} m^p}{2 M^p} \,.
\end{equation}
Thus partitions $\mathcal{N}$ with small number $r(\mathcal{N})$
make a small contribution in $T_p (Q)$.
At the second part of the proof we consider partitions $\mathcal{N}$
with large number of the sets $N_j$.

Let us estimate the sum $\sigma_2$.
Consider the sets $D_{i_1}, \dots, D_{i_{2p}}$.
Let $\mathcal{C}_j = N_j$.
So we get a partition of $[2p]$ onto the sets $\mathcal{C}_j$.
Using Lemma \ref{l:sophisticated}, we obtain
\begin{equation}\label{f:inter}
    \sigma_2
        \le
            \sum_{r=p-\lceil \d_0 \rceil}^{p}
                \sum_{\mathcal{N},\, r(\mathcal{N}) = r}\,
                    \sum_{\v{i} \in \mathcal{N}}
                        \left( \sum_{S^* \subseteq [2p],\, |S^*| = p}\,\, \per M_{\v{i}} (S^*) \right) \,.
\end{equation}
By Lemma \ref{l:sophisticated}
the summation in formula (\ref{f:inter}) is taken over
all sets $S^*$ such that $S^*$ contains an element from any set $N_j$.
Further let
$M_{\v{i}} (S^*) = (m_{\a\beta})$ be a matrix of size $p\m p$,
$m_{\a\beta} = |D_{i_\a} \bigcap D_{i_\beta}|$, $\a\in S^*$, $\beta\in \bar{S}^*$.
Let $M'_{\v{i}}$ be a matrix of the size $r\m 2p$,
$M'_{\v{i}} = (|D_\a \bigcap D_{i_\beta}|)_{\a\in S_{\v{i}}, \beta\in [2p]}$.
Using formula (\ref{f:permanent}), we get
\begin{equation}\label{f:X-}
    \per M_{\v{i}} (S^*)
        \le
            \prod_{\a \in S^*,\, i_\a \notin S_{\v{i}}}
                \left( \sum_{\beta \in \bar{S}^*} |D_{i_\a} \bigcap D_{i_\beta}| \right) \cdot \per M'_{\v{i}} \,.
\end{equation}
Applying the last inequality and a trivial estimate $|D_\l|\le s_2$, $\l \in \L_1$, we have
\begin{equation}\label{f:X}
    \per M_{\v{i}} (S^*)
        \le
            \prod_{\a \in S^*,\, i_\a \notin S_{\v{i}}}
                \left( \sum_{x\in \L_2} D_{i_\a} (x) \sum_{\beta \in \bar{S}^*} D_{i_\beta} (x) \right) \cdot \per M'_{\v{i}}
                    \le
                        p^{p-r} s_2^{p-r} \per M'_{\v{i}} \,.
\end{equation}
Using (\ref{f:permanent}), it is easy to see that
$$
    \per M'_{\v{i}}
        \le
            \pi (t_1,\dots,t_r) \prod_{\a\in S_{\v{i}}} \left( \sum_{\beta \in S_{\v{i}}} |D_\a \bigcap D_\beta| \right) \,,
$$
where the quantity  $\pi (t_1,\dots,t_r)$ was defined in Lemma \ref{l:columns}.
Using the bound for $\pi (t_1,\dots,t_r)$ from the lemma
and inequalities (\ref{f:inter}), (\ref{f:X}), we get
$$
    \sigma_2
        \le
            2^{5p} \max\{\, \d_0^{4\d_0}, 1 \,\}
                \cdot
                \sum_{r=p-\lceil \d_0 \rceil}^{p}
                    p^{p-r} s_2^{p-r}
                    \sum_{\mathcal{N},\, r(\mathcal{N}) = r}\,
                    \sum_{\v{i} \in \mathcal{N}}
                        \prod_{\a\in S_{\v{i}}} \left( \sum_{\beta \in S_{\v{i}}} |D_\a \bigcap D_\beta| \right) \,.
$$
If we make a permutation of components of the vector $\v{i}$ by
different parts $N_j$ of our partition $\mathcal{N}$ then
the set $S_{\v{i}}$ does not change.
Besides, if the lengths of sets $N_j$ are fixed then the set $S_{\v{i}}$
does not change after any permutation of these sets.
Hence
$$
    \sigma_2
        \le
            2^{5p} \max\{\, \d_0^{4\d_0}, 1 \,\}
                \cdot
                \sum_{r=p-\lceil \d_0 \rceil}^{p}
                    p^{p-r} s_2^{p-r}
                        \sum_{t_1+\dots+t_r = 2p} \frac{(2p)!}{t_1! \dots t_r!}\, \frac{1}{r!}\, r!
                        \cdot
$$
$$
                        \cdot
                        \left( \sum_{S\subseteq [s_1],\, |S| = r}\, \prod_{\a\in S}
                            \left( \sum_{\beta \in S} |D_\a \bigcap D_\beta| \right) \right)
                                \le
    2^{5p} \max\{\, \d_0^{4\d_0}, 1 \,\}
        p^{3p} s_2^p
        \cdot
$$
\begin{equation}\label{f:final2}
        \cdot
            \sum_{r=p-\lceil \d_0 \rceil}^{p}
                \left( \frac{1}{ps_2} \right)^r
                    \cdot
                        \left( \sum_{S\subseteq [s_1],\, |S| = r}\, \prod_{\a\in S}
                            \left( \sum_{\beta \in S} |D_\a \bigcap D_\beta| \right) \right) \,.
\end{equation}
Combining inequalities (\ref{f:final1}), (\ref{f:final2})
and formula (\ref{f:III}), we get inequality (\ref{e:T_p_intesect}).
This completes the proof.

To prove Theorem \ref{t:inverse2}
we need in a combinatorial lemma
and
a well--known lemma of E. Bombieri (see e.g. \cite{Vaughan}).

Let $p$ be a positive integer, and $A_1, \dots, A_p$
be a sequence of sets such that any two of them $A_i$ and $A_j$
either disjoint or equals.
By $\rho$ denote the number of different sets among $A_1,\dots, A_p$.
Let the set $A^*_1$ appears in the sequence $A_1,\dots, A_p$ exactly $l_1$ times,
$A^*_2$ --- exactly $l_2$ times, $\dots$, $A^*_\rho$ exactly $l_\rho$ times.

\Lemma
{\it
    Let $w$ be a positive integer, $2\le p \le a$,
    $\zeta \in (0,1]$ be a real number, and
    $S_1, \dots, S_q$ be some different sets, $|S_i| = p$, $S_i = \{ s^{(1)}_i, \dots, s^{(p)}_i \}$,
    $i=1,\dots,q$.
    Let also for all $i\in [q]$ and for all sets $S_i$, we have
    $s^{(j)}_i \in A_j$, $j=1,\dots,p$.
    Suppose that
    $$
        q \ge 2 \sum_{\o=\lceil \zeta p\rceil}^p\,
                \frac{(pw)^\o}{\o!} \sum_{n_1+\dots+n_\rho = p -\o,\, n_i \le l_i}
                    \frac{|A^*_1|^{n_1} \dots |A^*_\rho|^{n_\rho} }{n_1! \dots n_\rho!} \,.
    $$
    Then there are sets $S_{n_1}, \dots, S_{n_w}$ from the sequence $S_1, \dots, S_q$
    such that for an arbitrary $l=2,\dots,w$, we have
    $| (\bigcup_{i=1}^{l-1} S_{n_i} ) \bigcap S_{n_l} | \le \zeta p$.
}
\label{l:greed}
\\
\Proof
    We use a greedy algorithm.
    Let $S_{n_1} = S_1$.
    Suppose that sets $S_{n_1}, \dots, S_{n_{l-1}}$ have been constructed
    and find $S_{n_l}$.
    Let $C_l = \bigcup_{i=1}^{l-1} S_{n_i}$.
    Clearly, $|C_l| \le wp$.
    Let $C_l = C^*_1 \bigsqcup \dots \bigsqcup C^*_\rho$,
    where $C^*_i \subseteq A^*_i$, $i=1,\dots,\rho$.
    Let also $a_i = |A^*_i|$, $c_i = |C^*_i|$, $i=1,\dots,\rho$.
    The number of sets $E$ belong to $A^*_1 \bigsqcup \dots \bigsqcup A^*_\rho$,
    $|E| = p$ and such that $|E \bigcap C_l| \ge \zeta p$ does not exceed
$$
    \sigma
        :=
            \sum_{\o=\lceil \zeta p\rceil}^p\,
                \sum_{m_1+\dots+m_\rho = \o,\, m_i\le c_i }\,
                \sum_{n_1+\dots+n_\rho = p-\o,\, n_i\le \min\{ a_i-c_i,l_i \} }
                \binom{c_1}{m_1} \dots \binom{c_\rho}{m_\rho}
                \m
$$
$$
                \m
                \binom{a_1-c_1}{n_1} \dots \binom{a_\rho-c_\rho}{n_\rho}
            \le
                \sum_{\o=\lceil \zeta p\rceil}^p\,
                \sum_{m_1+\dots+m_\rho = \o}\,
                \sum_{n_1+\dots+n_\rho = p-\o,\, n_i \le l_i}
                    \frac{c_1^{m_1}\dots c_{\rho}^{m_\rho}}{m_1!\dots m_\rho!}
                    \cdot
                    \frac{a_1^{n_1}\dots a_{\rho}^{n_\rho}}{n_1!\dots n_\rho!}
                        \le
$$
$$
                        \le
                            \sum_{\o=\lceil \zeta p\rceil}^p\, \frac{(c_1 + \dots + c_\rho)^\o}{\o!}
                                 \cdot
                                    \sum_{n_1+\dots+n_\rho = p-\o,\, n_i \le l_i}
                                    \frac{a_1^{n_1}\dots a_{\rho}^{n_\rho}}{n_1!\dots n_\rho!}
                                \le
$$
$$
                                \le
                                    \sum_{\o=\lceil \zeta p\rceil}^p\, \frac{(pw)^\o}{\o!}
                                        \cdot
                                            \sum_{n_1+\dots+n_\rho = p-\o,\, n_i \le l_i}
                                            \frac{a_1^{n_1}\dots a_{\rho}^{n_\rho}}{n_1!\dots n_\rho!} = \sigma^* \,.
$$
    By assumption $q\ge 2 \sigma^*$.
    Hence $q\ge w$ and, consequently, $q-(l-1) > q - w \ge \sigma^*$.
    Thus there is a set $S_{n_l}$ from $S_1, \dots, S_q$ such that
    $S_{n_l}$ does not equal $S_{n_1}, \dots, S_{n_{l-1}}$ and such that
    $| (\bigcup_{i=1}^{l-1} S_{n_i} ) \bigcap S_{n_l} | \le \zeta p$.
    This completes the proof.

\Lemma {\bf (Bombieri)}
{\it
    Let $q$ be a positive integer, $\lambda >0 $ be a real number,
    $B$ be a finite set.
    Suppose that $B_1, \dots, B_q$ are subsets of $B$
    such that $|B_i| \ge \lambda |B|$.
    Then for all $t \le \lambda q$ there are different positive integers
    $j_1, \dots, j_t \in [q]$ such that
    $$
        | B_{j_1} \bigcap \dots \bigcap B_{j_t} | \ge \left( \lambda - \frac{t}{q} \right) \binom{q}{t}^{-1} |B| \,.
    $$
}
\label{l:Bombieri}

\Th
{\it
    Let $K,\eta>0$ be real numbers, $\eta \in (0,1/2]$, $p$ be a positive integer,
    and
    $\L \subseteq \mathbf{F}^n_2$ be an arbitrary set from the family $\mathbf{\L} (4p)$.
    Let also $Q$ be a subset of  $\L \dotplus \L$, $K^* := \max\{ 1,K \}$,
    $p \ge 2^{30} K^* / \eta$, and
    \begin{equation}\label{f:T_p_large}
        T_p (Q) \ge \frac{p^{2p} |Q|^p}{K^p} \,.
    \end{equation}
    Suppose that $p\le \log |\L| / \log \log |\L|$ and
    $$
        |Q|
            \ge
                2^{60+\frac{2}{\eta}} (K^*)^{17} p^3 |\L|
                    \cdot
                        \max \left\{ (2^{30} (K^*)^{11} p)^{\eta p} |\L|^{\eta} \log |\L|,
                            \exp \left(\frac{\log(2^{30} (K^*)^{20}) \log(\frac{p \log K^*}{\log p})}{\log(\frac{2^{-25} \eta p}{K^*})} \right) \right\}
    $$
    Then there are sets $\mathcal{L}_1, \mathcal{L}'_1, \dots, \mathcal{L}_h, \mathcal{L}'_h$
    from $\L$ such that
    $\mathcal{L}_i \bigcap \mathcal{L}'_j = \emptyset$,
    $\mathcal{L}_i + \mathcal{L}'_i \subseteq Q$, $i=1,\dots,h$, $j=1,\dots,h$,
    \begin{equation}\label{f:L_conditions}
        |\mathcal{L}_i| \ge \frac{\log (\frac{|Q|}{16(K^*)^9 |\L|})}{2^{10} \log (2^{20} K^*)},\, \quad
        |\mathcal{L}'_i| \ge \frac{1}{2^{10} p^2} \left( \frac{|Q|}{(K^*)^9 |\L|} \right)^{\eta}\,,
    \end{equation}
    $(\mathcal{L}_i + \mathcal{L}'_i) \bigcap (\mathcal{L}_j + \mathcal{L}'_j) = \emptyset$,
    $i,j=1,\dots,h$, $i\neq j$
    and
    \begin{equation}\label{f:large_part}
        \left| Q\bigcap \left( (\mathcal{L}_1 + \mathcal{L}'_1) \bigsqcup \dots \bigsqcup (\mathcal{L}_h + \mathcal{L}'_h) \right) \right|
            \ge
                \frac{|Q|}{16 (K^*)^9} \,.
    \end{equation}
    If $p$ is an arbitrary and
    \begin{equation}\label{f:m_cond2}
        \log \left( \frac{|Q|}{16(K^*)^9 p |\L|} \right) \ge 2^{20} \log (2^{10} K^*) \log p \,,
    \end{equation}
    then there are sets $\mathcal{L}_1, \mathcal{L}'_1, \dots, \mathcal{L}_h, \mathcal{L}'_h$
    from $\L$ satisfying (\ref{f:large_part}) and such that
    \begin{equation}\label{f:L_conditions2}
        |\mathcal{L}_i|
            \ge
                \min\{ 2^{-18} \frac{p}{K^*}, 2^{-5} \log  \left( \frac{|Q|}{16(K^*)^9 p} \right) \},\, \quad
        |\mathcal{L}'_i| \ge \frac{1}{32 p^2} \left( \frac{|Q|}{(K^*)^9 |\L|} \right)^{1/2}\,.
    \end{equation}
}
\label{t:inverse}

\Note
If $K = O(1)$ e.g. $K\le 1$ then inequalities (\ref{f:large_part}), (\ref{f:L_conditions2})
hold if we have more weaker bound than (\ref{f:m_cond2}), namely
$|Q| \ge  2^{60+\frac{2}{\eta}} (K^*)^{17} p^3 |\L|$.
\label{t:inverse2}
\\
{\bf Proof of the theorem.}
Let $m=|Q|$, $\beta_1 = 1/4$, $\beta_2 = 1/2$.
Let also
$$
    \mathbf{M}
        =
            2^{52+\frac{2}{\eta}} (K^*)^{17} p^3 |\L|
                    \cdot
                        \max \left\{ (2^{27} (K^*)^{11} p)^{\eta p} |\L|^{\eta} \log |\L|,
                            \exp \left(\frac{\log(2^{24} (K^*)^{20}) \log(\frac{p \log K^*}{\log p})}{\log(\frac{2^{-22} \eta p}{K^*})} \right) \right\} \,.
$$
Since $T_p (Q) \ge p^{2p} |Q|^p / K^p$, it follows that $D_p (Q) \ge p \log (p/K)$.
Using Theorem \ref{t:connected}  with parameters $d=2$ and $C=1/8$, we get
$(\beta_1, \beta_2)$--connected set $Q_1 \subseteq Q$ of degree $p$
such that
$m_1 := |Q_1| \ge m / (2K^9)$ and $T_p (Q_1) \ge p^{2p} m_1^p / K^p$.
Let $a=\lceil |\L|/2 \rceil$.
We have
\begin{equation}\label{f:half}
    \sum_{\t{\L} \subseteq \L,\, |\t{\L}| = a}\,
        \sum_{\l_1 \in \t{\L},\, \l_2 \in \L\setminus \t{\L}} Q_1 (\l_1 + \l_2)
            =
                2 \binom{|\L|-2}{a-1} |Q_1| \,.
\end{equation}
Using (\ref{f:half}), it is easy to see that there is a set
$\t{\L} \subseteq \L$, $|\t{\L}| = a$ such that
$
    | Q_1 \bigcap (\t{\L} + (\L\setminus \t{\L})) |
        \ge 2m_1 \binom{|\L|-2}{a-1} \binom{|\L|}{a}^{-1}
            = 2m_1 \frac{a (|\L|-a)}{|\L| (|\L|-1)} \ge m_1 / 2
$.
Put $\L_1 = \t{\L}$, $\L_2 = \L \setminus \t{\L}$ and
$Q_2 = Q_1 \bigcap (\L_1 + \L_2)$.
Certainly, we can find a set $Q_3 \subseteq Q_2$ such that $Q_3 = \lceil m_1 / 2\rceil$.
Let $m_3 = |Q_3|$.
Since the set $Q_1$ is
$(\beta_1, \beta_2)$--connected of degree $p$ and $C=1/8$, it follows that
$$
    T_p (Q_3) \ge 2^{-6p} \left( \frac{m_3}{m_1} \right)^{2p} T_p (Q_1) \ge \frac{p^{2p} m_3^p}{(2^7 K)^p} \,.
$$
Using notation of Proposition \ref{st:inverse2}, taking $M=2^{7} K$
and applying this Proposition to the set $Q_3 \subseteq \L_1 + \L_2$, we get
\begin{equation}\label{f:It}
    2^{5p} X p^{3p} s_2^p
                \cdot
                \sum_{r=p-\lceil \d_0 \rceil}^{p}
                \left( \frac{1}{ps_2} \right)^r
                    \cdot
                        \left( \sum_{S\subseteq [s_1],\, |S| = r}\, \prod_{\a\in S}
                            \left( \sum_{\beta \in S} |D_\a \bigcap D_\beta| \right) \right)
                                \ge
                                    \frac{p^{2p} m_3^p}{2(2^7 K)^p} \,.
\end{equation}
Recall that
the quantity $\d_0$ equals $\d_0 = \max\{ (p \log (2 e M) ) / \log (|Q_3|/(s_2 p)) , 1 \}$,
and the number $X$ is $\max\{\, \d_0^{4\d_0}, 1 \,\}$.
If $m\ge \mathbf{M}$ or $m$ satisfy (\ref{f:m_cond2}) then
$\d_0 \le \max \{ (p \log (2^{10} K) ) / (2\log p) , 1 \} \le p/2$
and $X^{1/p} \le 2^8 K$.
Let $K_1 = 2^{13} K X^{1/p} \le 2^{21} K^2$.
Suppose that either $m\ge \mathbf{M}$ or $m$ satisfy (\ref{f:m_cond2}).
Then
$m_3 \ge 2 K_1 p |\L|$.
Using the last inequality and (\ref{f:It}), we obtain
that there is a positive integer $p_1 \in [p-\lceil \d_0 \rceil, p]$
such that
\begin{equation}\label{f:tetrad'}
    \sum_{S\subseteq [s_1],\, |S| = p_1}\, \prod_{\a\in S}
         \left( \sum_{\beta \in S} |D_\a \bigcap D_\beta| \right)
            \ge
                \frac{m_3^{p_1}}{K_1^{p_1}} \,.
\end{equation}
Let $S\subseteq [s_1]$ be a set, $|S| = p_1$,
and $\a \in S$ be an arbitrary element of the set $S$.
Let also $\eps = 1/(16K_1)$.
If $M\le 1/2$ then $X = 1$,
and by inequality $p \ge 2^{30} K^* / \eta$, we get $\eps \ge 1/p_1$.
Suppose that $M>1/2$
and either $m\ge \mathbf{M}$ or $m$ satisfy (\ref{f:m_cond2}).
In the case the inequality $\eps \ge 1/p_1$
can be derived from the condition $p \ge 2^{30} K^* / \eta$.
A slightly more accurate computations show that in the both cases, we have $\eps \ge 16/(\eta p)$.
Define the sets
$$
    G_{S,\a} = \{~ x\in D_\a ~:~ \sum_{\beta \in S} D_\beta (x) \ge \eps p_1 ~\} \,.
$$
In other words, $G_{S,\a}$ is the set of $x$ from $D_\a$ such that
$x$ belongs to at least $\eps p_1$ the sets $D_\beta$, $\beta \in S$.
We have
$$
    \sum_{\beta \in S} |D_\a \bigcap D_\beta|
        =
            \sum_{x\in \L_2} D_\a (x) \sum_{\beta \in S} D_\beta (x)
                =
$$
\begin{equation}\label{f:sum_beta}
                =
                    \sum_{x\in G_{S,\a} } D_\a (x) \sum_{\beta \in S} D_\beta (x)
                        +
                            \sum_{x\notin G_{S,\a} } D_\a (x) \sum_{\beta \in S} D_\beta (x)
                                \le
                                    p_1 |G_{S,\a}| + \eps p_1 |D_\a| \,.
\end{equation}
Let $\mathcal{S}$ be the family of sets $S$,
$S\subseteq [s_1]$, $|S|=p_1$ such that
for any $S\in \mathcal{S}$ there is $\a \in S$
such that $|G_{S,\a}| \ge \eps |D_\a|$ and $|D_\a| \ge \eps m_3/s_2$.
Let also $\bar{\mathcal{S}}$ be the family of sets from $S$,
$S\subseteq [s_1]$, $|S|=p_1$ do not belong the family $\mathcal{S}$.
Let us prove that
\begin{equation}\label{f:1_est}
    \sigma_1
        :=
            \sum_{S \in \bar{\mathcal{S}}}\, \prod_{\a\in S}
                 \left( \sum_{\beta \in S} |D_\a \bigcap D_\beta| \right)
                    \le
                        \frac{m_3^{p_1}}{2K_1^{p_1}} \,.
\end{equation}
Let $Y(S) = \{ \a \in S ~:~ |G_{S,\a}| < \eps |D_\a| \}$, and $\bar{Y} (S) = S \setminus Y(S)$.
Using (\ref{f:sum_beta}), we get
$$
    \sigma_1
        =
            \sum_{S \in \bar{\mathcal{S}}}\, \prod_{\a\in \bar{Y} (S),\, |D_\a| < \eps m_3/s_2}
                \left( \sum_{\beta \in S} |D_\a \bigcap D_\beta| \right)
                    \cdot
                \prod_{\a\in Y (S)}
                \left( \sum_{\beta \in S} |D_\a \bigcap D_\beta| \right)
                    \le
$$
$$
                    \le
                        \sum_{l=0}^{p_1}
                            \sum_{S\subseteq [s_1],\, |S|=p_1,\, |Y(S)| = p_1 - l}
                                \left( \frac{\eps p_1 m_3}{s_2} \right)^{|\bar{Y} (S)|}
                                    \cdot
                                        \prod_{\a\in Y (S)} (2\eps p_1 |D_\a|)
    \le
$$
$$
    \le
        \sum_{l=0}^{p_1}
            \left( \frac{\eps p_1 m_3}{s_2} \right)^l
                (2\eps p_1)^{p_1-l} \binom{s_1 - (p_1 - l)}{l}
                    \sum_{S'\subseteq [s_1],\, |S'| = p_1 - l}\,
                        \prod_{\a \in S'} |D_\a|
                            \le
$$
$$
                            \le
                                2^{p_1} \eps^{p_1}
                                    \sum_{l=0}^{p_1} \left( \frac{\eps p_1 m_3}{s_2} \right)^l
                                        p_1^{p_1-l} \frac{s_1^l}{l!} \eps^{-l} \frac{1}{(p_1-l)!} m_3^{p_1-l}
    \le
        2 (2e)^{p_1} \eps^{p_1} m_3^{p_1}
            \sum_{l=0}^{p_1} \frac{p_1!}{l!(p_1-l)!}    
            =
$$
$$
            =
                2 (4e)^{p_1} \eps^{p_1} m_3^{p_1}       
                    \le
                        2^{4p_1 - 1} \eps^{p_1} m_3^{p_1}
                            =
                                \frac{m_3^{p_1}}{2K_1^{p_1}}
$$
and inequality (\ref{f:1_est}) is proved.
Hence
\begin{equation}\label{f:common_final}
 \sigma_2
    =
        \sum_{S \in \mathcal{S}}\, \prod_{\a\in S}
                 \left( \sum_{\beta \in S} |D_\a \bigcap D_\beta| \right)
                    \ge
                        \frac{m_3^{p_1}}{2K_1^{p_1}} \,.
\end{equation}

Consider the case $p\le \log |\L| / \log \log |\L|$.
Let $u_0 = [\log s_2]$,
and $\L^{(j)} = \{ \a \in [s_1] ~:~ 2^{j-1} \le |D_\a| \le 2^j \}$, $j=1,\dots, u_0$.
Applying inequality $\sum_{\a \in \L_1} |D_\a| \le m$,
we derive that $|\L^{(j)}| \le 2m 2^{-j}$.
Let $(j_1, \dots, j_{p_1})$ be a tuple from $[u_0]^{p_1}$.
By $\rho = \rho (j_1, \dots, j_{p_1})$ denote the number of                                             
different $j_t$
in $(j_1, \dots, j_{p_1})$
and let an element $j^*_1$ appears in the tuple
exactly $l_1$ times, an element $j^*_2$ appears exactly $l_2$ times, $\dots$,
and an element $j^*_\rho$ appears exactly $l_\rho$ times,
and all elements $j^*_1, j^*_2, \dots, j^*_\rho$ are different.
We have
$$
    \sigma_2
        \le
            \sum_{S \in \mathcal{S}}\, p_1^{p_1} \prod_{\a\in S} |D_\a|
                =
                    \frac{p_1^{p_1}}{p_1!} \sum_{S \in \mathcal{S}}\,
                        \sum_{\a_1,\dots,\a_{p_1} \mbox{--- different}}\,
                        S(\a_1) \dots S(\a_{p_1}) |D_{\a_1}| \dots |D_{\a_{p_1}}|
                =
$$
$$
                =
                        p_1^{p_1}
                        \sum_{S \in \mathcal{S}}\,
                        \sum_{j_1,\dots,j_{p_1}=1}^{u_0} \frac{1}{l_1!\dots l_\rho!} \m
$$
\begin{equation}\label{f:1.1}
                        \m \sum_{\a_1\in \L^{(j_1)}, \dots, \a_{p_1} \in \L^{(j_{p_1})},\,
                                \a_1,\dots,\a_{p_1} \mbox{--- different}}\,
                        S(\a_1) \dots S(\a_{p_1}) |D_{\a_1}| \dots |D_{\a_{p_1}}| \,.
\end{equation}
Using formula (\ref{f:1.1}), we obtain that
there is a tuple $(j_1,\dots,j_{p_1})$ such that
$$
    \sum_{S \in \mathcal{S},\, S = \{ s^{(1)}, \dots, s^{(p_1)} \},\, s^{(j)}\in \L^{(j)} }\,
        \prod_{\a\in S} |D_\a|
            \ge
                \frac{l_1! \dots l_\rho!}{p_1^{p_1} u_0^{p_1}} \cdot \frac{m_3^{p_1}}{4 K_1^{p_1}} \,.
$$
By the definition of the sets $\L^{(j)}$, we get
$$
    q_0:= | \{ S\in \mathcal{S} ~:~ S = \{ s^{(1)}, \dots, s^{(p_1)} \},\, s^{(j)}\in \L^{(j)} \} |
        \ge
            \frac{l_1! \dots l_\rho!}{p_1^{p_1} u_0^{p_1} 2^{j_1+\dots + j_{p_1}}} \cdot \frac{m_3^{p_1}}{4 K_1^{p_1}} \,.
$$
Using Dirichlet's principle, we obtain that there is $\a \in [s_1]$ such that
$$
    q:= | \{ S\in \mathcal{S} ~:~ S = \{ s^{(1)}, \dots, s^{(p_1)} \},\, s^{(j)}\in \L^{(j)},\, \a \in S \} |
            \ge
                \frac{l_1! \dots l_\rho!}{p_1^{p_1} u_0^{p_1} 2^{j_1+\dots + j_{p_1}}} \cdot \frac{m_3^{p_1}}{4 s_1 K_1^{p_1}} \,.
$$
Let $A_i = \L^{(j_i)}$, $i=1,\dots, p_1$, and $A^*_i = \L^{(j^*_i)}$, $i=1,\dots, \rho$.
We want to apply Lemma \ref{l:greed}
to the sets $A_i$, $A^*_i$
with parameter $w=[\log (m_3/s_2) / (\eps^2 p_1 \log (2^6/\eps^2)) ]$.
Since $m\ge \mathbf{M}$ and $m_3 \ge m/(8K^9)$, it follows that
\begin{equation}\label{f:check}
    q
        \ge
            \frac{l_1! \dots l_\rho!}{p_1^{p_1} u_0^{p_1} 2^{j^*_1 l_1+\dots + j^*_{\rho} l_\rho} } \cdot \frac{m_3^{p_1}}{4 s_1 K_1^{p_1}}
                \ge
                    2 \sum_{\o=\lceil \zeta p_1\rceil}^{p_1}\,
                    \frac{(p_1w)^\o}{\o!} \sum_{n_1+\dots+n_\rho = p_1 -\o,\, n_i \le l_i}
                    \frac{|A^*_1|^{n_1} \dots |A^*_\rho|^{n_\rho} }{n_1! \dots n_\rho!} = 2 \sigma^* \,.
\end{equation}
Indeed, by assumption $m\ge \mathbf{M} \ge p_1 w s_2$.
Hence
$$
    2^{j^*_1 l_1+\dots + j^*_{\rho} l_\rho} \sigma^*
        \le
            2^{p_1} \sum_{\o=\lceil \zeta p_1\rceil}^{p_1}\,
                    \frac{(p_1w)^\o}{\o!} \sum_{n_1+\dots+n_\rho = p_1 -\o,\, n_i \le l_i}
                        \frac{s_2^{\o} m^{p_1 - \o}}{n_1! \dots n_\rho!}
                         \le
$$
\begin{equation}\label{f:proverka_1}
                         \le
                            2^{p_1} m^{p_1}
                                \sum_{\o=\lceil \zeta p_1\rceil}^{p_1}\,
                                \frac{(p_1w)^\o}{\o!} \left( \frac{s_2}{m} \right)^{\o} \frac{\rho^{p_1-\o}}{(p_1 - \o)!}
                                \le
                                2^{4p_1} \left( \frac{\rho}{p_1} \right)^{p_1 (1-\zeta)} w^{\zeta p_1} m^{p_1 (1-\zeta)} s_2^{\zeta p_1}
\end{equation}
(we used the identity $1/(\o! (p-\o)!) = \binom{p}{w} / p!$
in the last inequality).
To check (\ref{f:check}) we need to verify inequality
\begin{equation}\label{f:1.10}
    m_3
        \ge
            2^7 u_0 K_1 w^{\zeta} p_1 m^{1-\zeta} s_1^{1/p_1} s_2^{\zeta}
                    \ge
            32 \left( \frac{\rho^{p_1(1-\zeta)} p_1^{\zeta p_1}}{l_1! \dots l_\rho!} \right)^{1/p_1}
                u_0 K_1 w^{\zeta} m^{1-\zeta} s_1^{1/p_1} s_2^{\zeta} \,.
\end{equation}
But the last inequality easily follows from $\eps \ge 16/(\eta p)$, $m_3 \ge m/(8K^9)$ and
$$
    m
        \ge
            \mathbf{M}
                \ge
                    2^{44} (K^*)^4 p |\L|^{1+\eta} \log |\L| ( 2^{27} (K^*)^{11} p )^{\eta p}
                        \ge
                            |\L| w p \cdot (s^{1/p_1} \log |\L|\, 2^{27} (K^*)^{11} p^2)^{1/\zeta} \,.
$$

Applying Lemma \ref{l:greed} to
the sets
$A_i$, $A^*_i$, we get new sets
$S^*_1, \dots, S^*_w \in \mathcal{S}$ such that for all
$l=2,3,\dots,w$, we have
$| (\bigcup_{i=1}^{l-1} S^*_{i} ) \bigcap S^*_{l} | \le \zeta p_1$.
Note that $|G_{S^*_i,\a}| \ge \eps |D_\a|$, $i=1,\dots,w$ and $|D_\a| \ge \eps m_3 / s_2$.
Applying Lemma \ref{l:Bombieri} with parameter $t=[\eps w/2]$
to the sets $G_{S^*_1,\a}, \dots, G_{S^*_w,\a} \subseteq D_\a$,
we get a set of indices $i_1 < \dots < i_t$ from $[w]$ such that
if $G^* = G_{S^*_{i_1},\a} \cap \dots \cap G_{S^*_{i_t},\a}$ then $|G^*| \ge \eps \binom{w}{t}^{-1} |D_\a|/2$.
Let $x$ be an arbitrary element of $D_\a$, and $\Gamma_i (x) = \{ \beta \in S^*_i ~:~ x\in D_\beta \}$.
Clearly, for any $x\in G_{S^*_{i},\a}$, we have $|\Gamma_i (x)| \ge \eps p_1$.
Let $E = \bigcup_{i=1}^w S^*_i$ and $I=\{ i_1, \dots, i_t \}$.
Obviously, $|E| \le w p_1$.
Consider the set
$$
    Z = \{~ x\in D_\a ~:~
            x \mbox{ belongs to at least } \frac{\eps p_1 t}{2} \mbox{ {\it different }sets } D_\beta,\, \beta \in E ~\} \,.
$$
Let us prove that $G^* \subseteq Z$.
Let $x\in G^*$.
Then $x$ belongs to the sets $D_\beta$, $\beta \in \bigcup_{i\in I} \Gamma_i (x)$.
Let us estimate the cardinality of $\bigcup_{i\in I} \Gamma_i (x)$.
We have
$$
    | \bigcup_{i\in I} \Gamma_i (x) |
        =
            | \bigcup_{i\in I\setminus \{ i_t \}} \Gamma_i (x) | + | \Gamma_{i_t} (x) |
                - | (\bigcup_{i\in I\setminus \{ i_t \}} \Gamma_i (x)) \bigcap \Gamma_{i_t} (x)|
                    \ge
$$
$$
                    \ge
                        | \bigcup_{i\in I\setminus \{ i_t \}} \Gamma_i (x) | + \eps p_1
                            - | (\bigcup_{i\in I\setminus \{ i_t \}} S^*_i ) \bigcap S^*_{i_t} |
                                \ge
                                    | \bigcup_{i\in I\setminus \{ i_t \}} \Gamma_i (x) | + \eps p_1 - \zeta p_1
                                        \ge \dots \ge \frac{\eps p_1 t}{2} \,.
$$
Whence $G^* \subseteq Z$ and, consequently,
$|Z| \ge |G^*| \ge \eps \binom{w}{t}^{-1} |D_\a|/2$.
Let $l = [\eps p_1 t/4]$.
Then
$$
    Z \subseteq \bigcup_{r_1,\dots,r_l \in E \mbox{ --- different }} \left( D_{r_1} \bigcap \dots \bigcap D_{r_l} \right) \,.
$$
Thus there is a tuple of indices $r_1< \dots < r_l$ from $E$ such that
$$
    |D_{r_1} \bigcap \dots \bigcap D_{r_l}|
        \ge
            \binom{p_1w}{l}^{-1} |Z|
                \ge
                    \binom{p_1w}{l}^{-1} \binom{w}{t}^{-1} \frac{\eps^2 m_3}{2s_2}
                        \ge
                            \frac{1}{2^{8} p^2} \left( \frac{m}{(K^*)^9 |\L|} \right)^{\eta} \,.
$$

Put $\mathcal{L}_1 = \{ \l_{r_1}, \dots, \l_{r_l} \}$ and
$\mathcal{L}'_1 = D_{r_1} \bigcap \dots \bigcap D_{r_l}$.
Then $\mathcal{L}_1 \bigcap \mathcal{L}'_1 = \emptyset$,
$\mathcal{L}_1 + \mathcal{L}'_1 \subseteq Q_3 \subseteq Q$
and
$$
  |\mathcal{L}_1| = l \ge \frac{\log (\frac{m_3}{s_2})}{32 \log (\frac{2^8}{\eps^2})}
    \ge \frac{\log (\frac{m}{8(K^*)^9 |\L|})}{2^{10} \log (2^{20} K^*)} \,.
$$
Now we can use an iterative procedure.
If $|\mathcal{L}_1 + \mathcal{L}'_1| \ge |Q_3| /2$ then we finish our procedure.
Otherwise consider the set $Q_3^{'} = Q_3 \setminus (\mathcal{L}_1 + \mathcal{L}'_1)$
and use our previous arguments.
We find sets
$\mathcal{L}_2 \subseteq \L_1$, $\mathcal{L}'_2 \subseteq \L_2$
such that
$\mathcal{L}_2 \bigcap \mathcal{L}'_2 = \emptyset$,
$\mathcal{L}_2 + \mathcal{L}'_2 \subseteq Q_3^{'} \subseteq Q$
and such that
$$
    |\mathcal{L}_2| \ge \frac{\log (\frac{|Q|}{16(K^*)^9 |\L|})}{2^8 \log (2^{12} K^*)},\, \quad
    |\mathcal{L}'_2| \ge \frac{1}{2^{10} p^2} \left( \frac{|Q|}{(K^*)^9 |\L|} \right)^{\eta}\,.
$$
By dissociativity of the sets $\L$ and $\L_1 \bigcap \L_2 = \emptyset$,
we get
$(\mathcal{L}_1 + \mathcal{L}'_1) \bigcap (\mathcal{L}_2 + \mathcal{L}'_2) = \emptyset$.
If $|\mathcal{L}_1 + \mathcal{L}'_1| + |\mathcal{L}_2 + \mathcal{L}'_2| \ge |Q_3| /2$
then we finish our algorithm.
At the end we construct sets
$\mathcal{L}_1, \mathcal{L}'_1, \dots, \mathcal{L}_h, \mathcal{L}'_h$
such that inequality (\ref{f:large_part}) holds.

We need to consider the case when (\ref{f:m_cond2}) holds
but either estimete $p\le \log |\L| / \log \log |\L|$ is not true or $m< \mathcal{M}$.
If inequality (\ref{f:m_cond2}) holds then $X=1$.
Using (\ref{f:common_final}) and a simple bound $|D_\a| \le s_2$, we obtain
that  the number of sets in the family $\mathcal{S}$ is at least
$m_3^{p_1}/(2 K_1^{p_1} p_1^{p_1} s_2^{p_1})$.
Applying condition (\ref{f:m_cond2}), we see that the last quantity is at least $1$.
Hence there is a set $S$ and a number $\a \in [s_1]$ such that
$|G_{S,\a}| \ge \eps |D_\a|$ and $|D_\a| \ge \eps m_3 / s_2$.
Let
$$
    Z = \{~ x\in D_\a ~:~
            x \mbox{ belongs to at least } \eps p_1 \mbox{ {\it different }sets } D_\beta,\, \beta \in S ~\} \,.
$$
Then $G_{S,\a} \subseteq Z$.
We have $[\eps p_1/2] \ge 2^{-18} \frac{p}{K^*} \ge 1$.
Put $l=\min \{ [\eps p_1/2], [\log (m_3/s_2) / 8] \}$.
We have
$$
    Z \subseteq \bigcup_{r_1,\dots,r_l \in S \mbox{ --- different }} \left( D_{r_1} \bigcap \dots \bigcap D_{r_l} \right) \,.
$$
Whence there is a tuple of indices $r_1< \dots < r_l$ from $S$ such that
$$
    |D_{r_1} \bigcap \dots \bigcap D_{r_l}|
        \ge
            \binom{\lceil \eps p_1 \rceil}{l}^{-1} |G_{S,\a}|
                \ge
                    \frac{\eps^2}{16^l} \frac{m_3}{s_2}
                        \ge
                            \frac{1}{16 p^2} \left( \frac{m}{(K^*)^9 |\L|} \right)^{1/2} \,.
$$
Put $\mathcal{L}_1 = \{ \l_{r_1}, \dots, \l_{r_l} \}$ and
$\mathcal{L}'_1 = D_{r_1} \bigcap \dots \bigcap D_{r_l}$.
Using the arguments as above, we get the required result.
This concludes the proof.

\Note
It is easy to see that the bound for the cardinalities of $\mathcal{L}_i$ from inequality (\ref{f:L_conditions})
is best possible.
We give a scheme of the proof of the last statement.
Let us preserve all notations of Theorem \ref{t:inverse}.
Let $K> 1$ be a fixed constant,
$\L_1, \L_2 \subseteq \L$, $\L_1 \bigcap \L_2 = \emptyset$.
Let
$Q \subseteq \L_1+ \L_2$ be a set which we will describe later, and $m:=|Q|$.
Let also $|\L_1| := s$, $|\L_2| = [mK/s]$.
Suppose that sets $D_\a \subseteq \L_2$, $\a = 1,\dots,s$
are random sets.
It means that for any $\a \in [s]$ an arbitrary element from $\L_2$,
belongs to set $D_\a$ with probability $1/K$.
Clearly, with positive probability, we have $|D_\a| \approx m/s$, $\a=1,\dots,s$
$|D_\a \bigcap D_\beta| \approx m/(sK)$, $\a \neq \beta$, $\a,\beta = 1,\dots,s$,
and
$$
    T_p (Q) \gg p^{2p} \sum_{S\subseteq [s],\, |S|=p} \prod_{\a \in S} \left( \sum_{\beta\in S} |D_\a \bigcap D_\beta| \right)
        \gg
            \frac{p^{2p} m^p}{K^p} \,.
$$
Thus inequality (\ref{f:T_p_large}) holds.
Nevertheless if $\mathcal{L}_1 \subseteq \L_1$, $\mathcal{L}_2 \subseteq \L_2$, $|\mathcal{L}_2| = l > 0$,
$\mathcal{L}_1 + \mathcal{L}_2 \subseteq Q$ then
$|\mathcal{L}_2| \le |D_{\a_1} \bigcap \dots \bigcap D_{\a_l}| \ll m /(sK^l)$
and we get a bound $l\ll \log (m/s) / \log K$.



Theorem \ref{t:inverse} has a simple corollary.

\Pred
{\it
    Let $K,\eta>0$ be real numbers, $\eta \in (0,1/2]$, $K\ge 1$,
    $p,d$ be positive integers,
    $d\ge 3$,
    and
    $\L \subseteq \mathbf{F}^n_2$ be an arbitrary set, $\L \in \mathbf{\L} (2dp)$.
    Let also $Q$ be a subset of $d \dot{\L}$, $|\L| \ge 8 d^2$,
    $p \ge 2^{50+8d} K^d / \eta$ and
    \begin{equation}\label{conD:T_p_large_d}
        T_p (Q) \ge \frac{p^{dp} |Q|^p}{K^{(d-1)p}} \,.
    \end{equation}
    Suppose that $p\le \log |\L| / \log \log |\L|$ and
    $$
        |Q|
            \ge
                2^{60+50d+\frac{2}{\eta}} M^{17} K^{2d} p^3 d^{-d} |\L|^{d-1}
                    \m
    $$
    $$
                    \m
                        \max \left\{ (2^{30} M^{11} p)^{\eta p} |\L|^{\eta} \log |\L|,
                            \exp \left(\frac{\log(2^{30} M^{20}) \log(\frac{p \log M}{\log p})}{\log(\frac{2^{-25} \eta p}{M})} \right) \right\}\,,
    $$
    where $M = 2^{13} (8K)^{d-1}$.
    Then there are sets $\mathcal{L}, \mathcal{L}' \subseteq \L$
    and elements $\l_1+\dots+\l_{d-2}$ from $\L$ such that
    $\mathcal{L}_i \bigcap \mathcal{L}'_j = \emptyset$,
    $\l_i \notin \mathcal{L}, \mathcal{L}'$,
    \begin{equation}\label{f:L_conditions_d*}
        |\mathcal{L}| \ge \frac{\log (\frac{|Q|d^d}{2^{140+80d} K^{3d} |\L|^{d-1}})}{2^{10} \log (2^{40} 8^d K^d)},\, \quad
        |\mathcal{L}'| \ge \frac{1}{2^{10} p^2} \left( \frac{|Q| d^d}{2^{140+80d} K^{3d} |\L|^{d-1}} \right)^{\eta}\,,
    \end{equation}
    and
    \begin{equation}\label{f:t_2inverse_conclusion}
        \l_1 + \dots + \l_{d-2} + \mathcal{L} + \mathcal{L}' \subseteq Q \,.
    \end{equation}
}
\label{pr:2inverse}
\Proof
Let $m=|Q|$, $\beta_1 = 4^{-d}$, $\beta_2 = 4^{-d} + 1/\sqrt{m}$, and $a=\lceil |\L|/d \rceil$.
Since $T_p (Q) \ge p^{dp} |Q|^{p} / K^{(d-1)p}$, it follows that $D_p (Q) \ge (d-1) p \log (p/K)$.
Using Theorem \ref{t:connected}  with parameters $d$ and $C=2^{-6}$, we get
$(\beta_1, \beta_2)$--connected set $Q_1 \subseteq Q$ of degree $p$ such that
$m_1 := |Q_1| \ge m / (d K^{2(d-1)} )$ and $T_p (Q_1) \ge p^{dp} m_1^p / K^{(d-1)p}$.
Let also $a_i = a$, $i=1,\dots,d-1$ and $a_d = |\L| - \sum_{i=1}^{d-2} a_i$.
Since $|\L| \ge 8 d^2$, it follows that $|\L|/(2d) \le a_d \le |\L|/d$.
It is easy to see that
\begin{equation}\label{f:halF*}
    Q(x)
        =
            \left( \frac{d! (|\L| - d)!}{(a_1 - 1)! \dots (a_d - 1)!} \right)^{-1}
                \sum_{S_1,\dots,S_d,\, |S_i| = a_i,\, \bigsqcup_{i=1}^d S_i = \L}\,
                    \left( Q \bigcap (S_1 + \dots + S_d) \right) (x) \,.
\end{equation}
Using formula (\ref{f:halF*}), we obtain that there is a
tuple of disjoint sets $S_1,\dots, S_d \subseteq \L$ such that
$$
    | Q_1 \bigcap (S_1 + \dots + S_d) |
        \ge
            m_1 d! \frac{(|\L| - d)!}{(a_1 - 1)! \dots (a_d - 1)!} \left( \frac{|\L|!}{a_1! \dots a_d!} \right)^{-1}
                =
$$
\begin{equation}\label{}
                =
                    m_1 d! \frac{a_1 \dots a_d}{|\L| (|\L|-1) \dots (|\L|-d+1)}
                        \ge
                            \frac{1}{2} e^{-d} m_1 \,.
\end{equation}
Put $Q_2 = Q_1 \bigcap (S_1 + \dots + S_d)$.

Let $d_1$ be a positive integer, $d_1 \le d$,
$l_1,\dots, l_{d_1}$ be different numbers from $[d]$,
$L = \{ l_1,\dots, l_{d_1} \}$, $\ov{L} = [d] \setminus L$.
Let also
$w_{l_i} \in S_{l_i}$ be arbitrary elements, $i\in [d_1]$,
$\v{w} = ( w_{l_1},\dots, w_{l_{d_1}} )$ be a vector, and
$W = \{ w_{l_1},\dots, w_{l_{d_1}} \}$.
Define the sets $D(W)$, $Q(W)$
$$
    D(W) = \{ \sum_{i\in \ov{L}} \l_i ~:~  \sum_{i\in \ov{L}} \l_i + \sum_{i\in L} w_i \in Q' \},\, \quad
    Q(W) = \{ q\in Q' ~:~ q = \sum_{i\in \ov{L}} \l_i + \sum_{i\in L} w_i \}
$$
Clearly, $D(W) = Q(W) + \sum_{i\in L} w_i$.
We shall write $D(\v{w})$, $Q(\v{w})$ instead of $D(W)$, $Q(W)$.
By assumption the set $\L$ belongs to the family $\mathbf{\L} (2dp)$.
Using this, it is easy to see that the sets $Q(W \bigcup \{ l_1 \})$ and $Q(W \bigcup \{ l_2 \})$,
$\l_1\neq \l_2$ are disjoint.
Besides, $Q(W \bigcup \{ l \}) \subseteq Q(W)$.
Whence, for all $x\ge 1$, we have
\begin{equation}\label{f:2inverse_10}
    \sum_{\l} |Q(W \bigcup \{ \l \})|^x
        \le
            |Q(W)|^{x-1} \sum_\l |Q(W \bigcup \{ \l \})|
                = |Q(W)|^x \,.
\end{equation}
Let $x_1,x_2 \ge 1$ be arbitrary numbers.
Using bound (\ref{f:2inverse_10}) and Cauchy--Schwartz, we get
\begin{equation}\label{f:2inverse_11}
    \sum_{\l} |Q(W_1 \bigcup \{ \l \})|^{x_1/2} |Q(W_2 \bigcup \{ \l \})|^{x_2/2}
        \le
            |Q(W_1)|^{x_1/2} |Q(W_1)|^{x_2/2} \,.
\end{equation}
Clearly, there is an analog of formula (\ref{f:2inverse_11})
for larger number of sets $Q(W_i \bigcup \{ \l \})$.

By $Q'_2$ denote the union of the sets $Q_2 (\v{a})$, $\v{a} \in S_1 \m \dots \m S_{d-2}$
such that  $|Q_2 (\v{a})| \ge |Q_2| /(4|S_1|\dots |S_{d-2}|)$.
Then $|Q'_2| \ge |Q_2|/2$.
Certainly, we can find a set $Q' \subseteq Q'_2$ such that $|Q'| = \lceil 4^{-d} m_1 \rceil$
and such that $Q' = \bigsqcup \t{Q}_2 (\v{a})$, $\t{Q}_2 (\v{a}) \ge |Q_2| /(16|S_1|\dots |S_{d-2}|)$.
Let $m' = |Q'|$.
Since the set $Q_1$ is $(\beta_1, \beta_2)$--connected of degree $p$ and $C=2^{-6}$, it follows that
\begin{equation}\label{f:great_or_equal}
    T_p (Q')
        \ge
            2^{-12p} \left( \frac{m'}{m_1} \right)^{2p} T_p (Q_1)
                \ge
                    \frac{p^{dp} m'^p}{2^{12 p} 4^{dp} K^{(d-1)p}} \,.
\end{equation}

Consider the equation
\begin{equation}\label{f:2_inverse_(1)}
    q_1 + \dots + q_{2p} = 0 \,,
\end{equation}
where $q_i \in Q'$, $i=1,\dots, 2p$.
By $\sigma'$ denote the number of solutions
of
(\ref{f:2_inverse_(1)}).
Let $\v{a}_1,\dots, \v{a}_{d-2}$ be arbitrary vectors from $S_1\m \dots \m S_{d-2}$,
and let $\v{v} = ( \v{a}_1,\dots, \v{a}_{2p} )$.
Denote by $\sigma(\v{v}) = \sigma (\v{a}_1,\dots, \v{a}_{2p})$ the set of
solutions of equation (\ref{f:2_inverse_(1)}) such that  $q_i \in Q(\v{a}_i)$, $i\in [2p]$.
Further by $\mathcal{M}$ denote the family of partitions of the segment $[2p]$
onto $p$ sets $\{ C_1,\dots, C_p \}$, $|C_j| = 2$, $j\in [p]$.
Let also $\mathcal{V}$ be the collection of all partitions $\{ \mathcal{M}_1,\dots,\mathcal{M}_{d-2} \}$,
$\mathcal{M}_i \in \mathcal{M}$, $i\in [d-2]$.
Clearly, the total number of different tuples $\{ C_1,\dots, C_p \}$
in $\mathcal{V}$ does not exceed $p^{p(d-2)}$.
By definition a vector $\v{v} = (\v{a}_1,\dots, \v{a}_{2p})$
belongs to $\mathcal{V}$ if for any $j=1,\dots, d-2$ and for all
set $C$ of partition $\mathcal{M}_j$, $C = \{ \a,\beta \}$, we have $\l_{\a} = \l_{\beta}$.
Obviously
\begin{equation}\label{f:2inverse_tilde}
    \sigma'
        \le
            \sum_{\mathcal{V}} \sum_{(\v{a}_1,\dots, \v{a}_{d-2}) \in \mathcal{V}} |\sigma ( (\v{a}_1,\dots, \v{a}_{d-2}) )|
                =
            \sum_{\mathcal{V}} \sum_{\v{v} \in \mathcal{V}} |\sigma (\v{v})| \,.
\end{equation}
Using Lemma \ref{l:conv}, we get
\begin{equation}\label{f:2inverse_ttilde}
    |\sigma (\v{v})|
        \le
            \left( \prod_{i=1}^{2p} T_p (\v{a}_i) \right)^{1/2p} \,.
\end{equation}
Suppose that for all vectors
$\v{a}$ from $S_1\m \dots \m S_{d-2}$, we have
\begin{equation}\label{f:2inverse_I}
    T_p (\v{a}) \le \frac{p^{2p} |Q(\v{a})|^p}{M^p} \,,
\end{equation}
where $M = 2^{13} (8K)^{(d-1)}$.
By the last inequality and (\ref{f:2inverse_tilde}), (\ref{f:2inverse_ttilde}), we obtain
\begin{equation}\label{}
    \sigma'
        \le
            \frac{p^{2p}}{M^p}
                \sum_{\mathcal{V}} \sum_{(\v{a}_1,\dots, \v{a}_{d-2}) \in \mathcal{V}}
                    \prod_{i=1}^{2p} | Q(\v{a}_i) |^{1/2} \,.
\end{equation}
Using formulas (\ref{f:2inverse_10}), (\ref{f:2inverse_11}) several times, we get
$$
    \sigma'
        \le
            \frac{p^{2p}}{M^p}  \sum_{\mathcal{V}} m^p
                \le
                \frac{p^{dp} m^p}{M^p} \,.
$$
We obtain a contradiction with (\ref{f:great_or_equal}).
Hence there is vector $\v{a}$ from $S_1\m \dots \m S_{d-2}$
such that  inequality (\ref{f:2inverse_I}) does not hold
and
$$
    |Q(\v{a})|
        \ge
            \frac{|Q_2|}{4|S_1|\dots |S_{d-2}|}
                \ge
                    \frac{m}{64d e^d (4K)^{2(d-1)} |S_1|\dots |S_{d-2}| }
                        \ge
                            \frac{m d^d}{2^{50d} K^{2d} |\L|^{d-2}} \,.
$$
Let $\v{a} = (a_1,\dots, a_{d-2})$.
Put $\l_i = a_i$, $i\in [d-2]$.
Applying Theorem \ref{t:inverse} to the set $Q(\v{a}) \subseteq S_{d-1}+S_d$,
we get sets $\mathcal{L}$, $\mathcal{L}'$ such that
inequalitiy (\ref{f:L_conditions_d*}) holds
and inclusion (\ref{f:t_2inverse_conclusion}) is true.
This completes the proof.

\refstepcounter{section}
\label{appendix}

\begin{center}
{\large
{\bf \arabic{section}. Appendix.}}
\end{center}

In the section we prove an analog of Theorem \ref{t:main}
for an arbitrary Abelian group $G$.

\Th
 \label{t:main_G}
 {\it
    Let $\d,\a$ be real numbers, $0< \a \le \d$,
    $A$ be a subset of $G$, $|A| = \d |G|$,
    $k\ge 2$ be a positive integer, and the set $\r_\a$
    as in (\ref{f:R_def}).
    Suppose that $B\subseteq \r_\a$ be an arbitrary set.
    Then
    $$
        T_k (B) \ge \frac{\d \a^{2k}}{\d^{2k}} |B|^{2k} \,.
    $$
 }
\Proof
Let $r\in \F{G}$.
Define the quantity $\theta(r) \in \mathbf{S}^1$
by the formula $\F{A} (r) = |\F{A} (r)| \theta(r)$.
We have
$$
    \a N |B| \le \sum_{r\in B} |\F{A} (r)| = \sum_x \sum_r B(r) \theta^{-1} (r)  e(-r\cdot x) \,.
$$
Using H\"{o}lder's inequality, we obtain
\begin{equation}\label{f:Gelder_G}
    (\a N |B|)^{2k}
        \le
            \sum_x \left| \sum_r B(r) \theta^{-1} (r)  e(-r\cdot x) \right|^{2k} \cdot \left( \sum_x A(x) \right)^{2k-1}
                = N T_k (B\cdot \theta^{-1}) (\d N)^{2k-1} \,.
\end{equation}
Let us prove a simple lemma.

\Lemma
{\it
    Let $f: G \to \mathbb {C}$ be an arbitrary function,
    and $k\ge 2$ be a positive integer.
    Then $T_k (f) \le T_k (|f|)$.
}
\\
{\bf Proof of the lemma.}
Using the triangle inequality, we get for any functions
$g,h : G \to \mathbb {C}$ the following holds
$|(g*h) (x)| \le (|g|*|h|) (x)$, $x\in G$.
By definition of $T_k (f)$ we obtain the required result.

Using the last lemma and (\ref{f:Gelder_G}), we get
$T_k (B) \ge T_k (B\cdot \theta^{-1}) \ge \frac{\d \a^{2k}}{\d^{2k}} |B|^{2k}$.
This concludes the proof.

\end{document}